\documentclass[11pt,letterpaper]{article}
\usepackage[T1]{fontenc}
\usepackage{lmodern,microtype}
\usepackage[margin=1in]{geometry}
\usepackage{amsmath,amssymb,amsthm,mathtools,booktabs,array}
\usepackage{enumitem,needspace}

\usepackage{tikz}
\usetikzlibrary{arrows.meta,calc,decorations.markings,positioning}
\usepackage{caption}
\usepackage[section]{placeins}
\usepackage{flafter}
\definecolor{smallfill}{RGB}{242,214,208}
\definecolor{mediumfill}{RGB}{223,229,224}
\definecolor{largefill}{RGB}{207,224,239}
\definecolor{figureaccent}{RGB}{41,77,111}
\tikzset{
  edge/.style={draw=black!65,line width=.5pt},
  contact/.style={draw=black!48,line width=.35pt},
  strongedge/.style={draw=figureaccent,line width=1.3pt},
  disk/.style={draw=black!80,line width=.45pt},
  vertex/.style={circle,draw=black!80,line width=.5pt,fill=white,
    minimum size=5.8mm,inner sep=0pt,font=\small},
  typeS/.style={fill=smallfill},typeM/.style={fill=mediumfill},
  typeL/.style={fill=largefill},
  arrow/.style={-{Stealth[length=2mm]},line width=.65pt},
  twoway/.style={{Stealth[length=2mm]}-{Stealth[length=2mm]},line width=.65pt},
  annotation/.style={font=\footnotesize,align=center},
  panel/.style={font=\small,align=center},
  tick/.style={line width=.5pt}
}

\usepackage[colorlinks=true]{hyperref}
\hypersetup{pdftitle={Optimal Radius Ratio of Compact Nonuniform Packings},pdfauthor={Zhen Zhang},pdfsubject={Whole-plane optimality and uniqueness of packing 53}}
\setlist[enumerate]{label=\textup{(\roman*)},leftmargin=2.1em,itemsep=3pt,topsep=5pt}
\numberwithin{equation}{section}
\newtheorem{theorem}{Theorem}[section]
\newtheorem{proposition}[theorem]{Proposition}
\newtheorem{lemma}[theorem]{Lemma}
\newtheorem{corollary}[theorem]{Corollary}
\theoremstyle{definition}
\newtheorem{definition}[theorem]{Definition}
\newtheorem{problem}[theorem]{Problem}
\theoremstyle{remark}
\newtheorem{remark}[theorem]{Remark}
\newcommand{\R}{\mathbb R}
\newcommand{\PP}{\mathcal P}
\newcommand{\word}[1]{\mathsf{#1}}
\newcommand{\II}{\mathcal I}
\newcommand{\calA}{\mathcal A}
\title{Optimal Radius Ratio of Compact Nonuniform Packings}
\author{Zhen Zhang\\[0.35em]
{\small Yau Mathematical Sciences Center, Tsinghua University, Beijing, China}\\
{\small\href{mailto:zz628@cornell.edu}{\texttt{zz628@cornell.edu}}}}
\date{November 19, 2026}
\begin{document}
\maketitle
\begin{abstract}
We study the ratio of the smallest to the largest radius in a nonuniform disk packing whose tangency graph triangulates the plane. We give a computer-assisted proof that the largest possible ratio is $s=0.651050185882609\ldots$, attained uniquely up to similarity by packing~53 of Fernique, Hashemi, and Sizova. In fact, every nonuniform packing with radius ratio at least $0.65105$ is similar to this packing. No periodicity or restriction on the number of radii is assumed; the radii are required only to have a positive infimum and a finite supremum. Following the finite-subgraph exclusion program of Connelly and Zhang, we first restrict the degree-five flowers and their pairings. The remaining exceptional degree-seven corona forces seventeen complete flowers on an unfolded disk with thirty-seven labels. An exact interval certificate excludes their angle equations. A comparison of logarithmic radii on a recurrent network then forces the radii and the three local degree words. We identify the entire triangulation by reducing these words to the square--octagon tiling and showing that its decoration is forced.
\end{abstract}

\section{Introduction}\label{sec:introduction}
The search for dense packings of disks with nearly equal radii is a classical problem. In \emph{Regular Figures} (1964), L\'aszl\'o Fejes T\'oth~\cite[\S27]{ft1964} studied how disks with radii in a prescribed interval could be chosen and arranged to maximize packing density. A central question is how little variation in disk size is needed to improve on the hexagonal density $\pi/\sqrt{12}$. His constructions remained benchmarks for decades; see Fernique~\cite[Section~4]{fernique-survey} for the subsequent development of this question.

Triangulated packings form a natural class in this search. A \emph{triangulated disk packing}, also called a \emph{compact packing}, is a family of disks with disjoint interiors such that joining the centers of tangent disks gives a locally finite triangulation of the whole plane. Fejes T\'oth~\cite[Theorem~1]{ft1984} proved that, when the radii are bounded above and away from zero, its lower density\footnote{The lower density is the limit inferior of the covered area fraction in disks centered at a fixed point as their radius tends to infinity.} is at least $\pi/\sqrt{12}$. For a packing with unequal radii and two independent translation periods, the inequality is strict~\cite[Corollary~2.4]{connelly-pierre}. Thus doubly periodic triangulated examples with large radius ratio give packings of nearly equal disks that are denser than any packing of equal disks. This density advantage motivates the search for the largest possible radius ratio in the triangulated class.

Let $\PP=\{\overline B(x_v,r_v):v\in V\}$ be such a packing, where $x_v\in\R^2$ and $r_v>0$ are the center and radius of the disk indexed by $v$. Here $\overline B(x,r)$ denotes the closed Euclidean disk. We assume throughout that
\begin{equation}\label{eq:bounded-radii}
0<\inf_{v\in V}r_v\le\sup_{v\in V}r_v<\infty,
\end{equation}
and define its \emph{radius ratio}, or \emph{uniformity}, by
\begin{equation}\label{eq:ratio}
q(\PP)=\frac{\inf_v r_v}{\sup_v r_v}.
\end{equation}
We call the packing \emph{nonuniform} if its radii are not all equal. The infimum and supremum in \eqref{eq:ratio} need not be attained.

This paper determines the sharp radius ratio among all nonuniform compact packings satisfying \eqref{eq:bounded-radii}. We prove that packing~53 of Fernique, Hashemi, and Sizova~\cite{fhs} is the unique extremal packing up to similarity. The radii are allowed to take infinitely many distinct values, and no periodicity assumption is imposed. More strongly, every such packing with radius ratio at least $0.65105$ is similar to packing~53. Thus the radius bound forces exactly three radii, their proportions, the entire tangency graph, and two independent translation periods.

Among the examples in \emph{Regular Figures} is a triangulated two-size packing with radius ratio approximately $0.63756$. Kennedy~\cite{kennedy} proved that exactly nine ratios admit compact packings with two sizes and that this example attains the largest of them. For three sizes, Messerschmidt~\cite{messerschmidt-three} established finiteness of the admissible normalized radius triples, and Fernique, Hashemi, and Sizova~\cite{fhs} determined their exact number, namely $164$. Their packing~53 improves the radius ratio to approximately $0.65105$ and is shown in Figure~\ref{fig:reference-packing}. More generally, Messerschmidt~\cite{messerschmidt-finite} proved finiteness of the normalized radius tuples for every fixed number of sizes.

The connection with density has also led to results with the available radii prescribed. B\'edaride and Fernique~\cite{bedaride-fernique} gave a unified proof that a compact packing attains maximal density at each of the nine binary compact ratios, and Fernique and Pchelina~\cite{fernique-pchelina} studied the ternary case. Our theorem addresses the uniformity question by maximizing \eqref{eq:ratio} over the triangulated class, with the radii themselves allowed to vary.

Connelly and Zhang~\cite[Conjecture~1]{cz} formulated the specific assertion that packing~53 is optimal and proposed a finite-subgraph exclusion program for proving it. Theorem~\ref{thm:optimal} proves this conjecture and answers the triangulated part of the uniformity question in Fernique~\cite[Section~4, Open Question~4]{fernique-survey}. We follow this program: finite exclusions and a comparison of radii force the local words, and a final geometric argument identifies the whole packing.

Let $s$ be the unique root in $(0.65105018,0.65105019)$ of
\begin{equation}\label{eq:polynomial}
\begin{split}
P(X)={}&89X^8+1344X^7+4008X^6-464X^5-2410X^4\\
       &\hspace{18mm}+176X^3+296X^2-96X+1.
\end{split}
\end{equation}
Let $r\in(0.83420,0.83440)$ be the medium radius in packing~53, with the largest radius normalized to $1$. Thus
\begin{equation}\label{eq:reference-numbers}
s=0.651050185882609\ldots,
\qquad r=0.834306042853017\ldots.
\end{equation}
These numbers are fixed reference radii. Their exact angle characterization is given in Section~\ref{sec:reference}; the displayed decimal expansions are not used as exact inputs to the verification.

The uniqueness of a nonuniform compact packing for the fixed radii $s,r,1$ was already recorded by Fernique~\cite[p.~2]{fernique2019}, with attribution to~\cite{fhs}. The catalogue in~\cite[Supplementary material~1, Apps.~A--B]{fhs} places case~53 in its doubly periodic class. The additional assertion proved here is that the radius bound forces an arbitrary set of radii to be proportional to this fixed triple. We also give a direct geometric identification of the packing from its three local degree words in Section~\ref{subsec:identification}.

\input{figures/01_reference-packing.tikz}

\subsection{Main results}
Throughout the paper we use the rational cutoff
\begin{equation}\label{eq:threshold}
q_0:=\frac{65105}{100000}=0.65105.
\end{equation}
The local exclusions and the exceptional-corona exclusion are all established at this value. It satisfies
\begin{equation}\label{eq:threshold-order}
0.645<q_0<s.
\end{equation}

Two packings are \emph{similar} if one is the image of the other under a map $x\mapsto\lambda Ux+t$, where $\lambda>0$, $U$ is an orthogonal linear map of $\R^2$, and $t\in\R^2$. Thus reflections are allowed. A packing is \emph{doubly periodic} if it is invariant under translations by two linearly independent vectors.

\begin{theorem}\label{thm:optimal}
Among nonuniform triangulated disk packings of the plane satisfying \eqref{eq:bounded-radii}, the largest radius ratio is $s$. Equality holds exactly for the packings similar to packing~53 of Fernique, Hashemi, and Sizova.
\end{theorem}

We obtain the following stronger statement. The degree of a disk means the degree of its vertex in the tangency graph; write $d(v)$ for the degree of $v$. A \emph{degree word} is the cyclic sequence of the degrees of its neighbors, considered up to rotation and reversal. The letters $S,M,L$ will always mean degrees $5,6,7$, respectively; they do not initially prescribe radii.

\begin{theorem}\label{thm:structure}
Let $\PP$ be a nonuniform triangulated disk packing satisfying \eqref{eq:bounded-radii}. If $q(\PP)\ge q_0$, then $\PP$ is similar to packing~53. In particular, it is doubly periodic, every vertex has degree $5,6$, or $7$, and there is a number $\lambda>0$ such that
\begin{equation}\label{eq:forced-radii}
r_v=\lambda\begin{cases}
s,&d(v)=5,\\
r,&d(v)=6,\\
1,&d(v)=7,
\end{cases}
\end{equation}
All three degrees occur. The only degree words are
\begin{equation}\label{eq:ideal-words}
\begin{array}{c|c}
\text{degree of the center}&\text{neighbor word}\\\hline
5&\word{LLMLM}\\
6&\word{LMLSLS}\\
7&\word{SLSMSMM}.
\end{array}
\end{equation}
Consequently $q(\PP)=s$.
\end{theorem}

Theorem~\ref{thm:structure} gives a gap: a nonuniform packing cannot have radius ratio in $[q_0,s)$. It also identifies every packing in the entire range $q\ge q_0$. The three-radius and double-periodicity conclusions require neither a hypothesis on the number of radii nor a periodicity hypothesis.

\subsection{Method of proof}
The proof consists of local exclusions, a comparison of radii, and a global identification. The local part uses the degree restrictions of Connelly and Zhang~\cite[Propositions~7.2--7.5]{cz} together with two finite exclusions at $q_0$. These force every degree-five flower to have word $\word{LLMLM}$ and give a pairing across its unique consecutive $LL$ edge. They also exclude a degree-six vertex with three degree-five neighbors.

There are then two possibilities for a degree-seven vertex with three degree-five neighbors. One is the reference word $\word{SLSMSMM}$. The other, $\word{SMSMSMM}$, is exceptional. We show that this second word forces a finite system of seventeen flower equations on thirty-seven independent radius variables. A complete interval certificate proves that the system has no solution with all radii in $[q_0,1]$. The construction is an unfolded disk: its map into a putative packing need not be injective.

After this exclusion, assigning reference radii $s,r,1$ according to degree gives an angle sum at least $2\pi$ at every vertex. The difference between actual and reference logarithmic radii is therefore a bounded superharmonic function on a positive symmetric conductance network. The Euclidean geometry of the packing makes this network recurrent. Its bounded superharmonic functions are constant, which forces the radii and the equality words in \eqref{eq:ideal-words}.

This comparison is a bounded-radius instance of Dubejko's theorem for positive labels on recurrent triangulations~\cite[Theorem~3.1]{dubejko}. The connection between recurrence and the type of an infinite circle packing was developed by He and Schramm~\cite{he-schramm}. We include the comparison and recurrence arguments in the form needed here.

Finally, merging paired degree-five flowers and paired $LMM$ triangles gives a plane map with one quadrilateral and two octagons at every vertex. This is the square--octagon tiling. Its restoration to the original triangulation is determined by one diagonal choice, and the two initial choices are equivalent. The fixed edge lengths then determine the packing up to similarity. This last step uses the complete equality words obtained from the comparison; excluding the exceptional word alone does not yet give them.

Section~\ref{sec:definitions} gives the angle and comparison formulas. Section~\ref{sec:local} establishes the local degree restrictions, and Section~\ref{sec:exceptional} excludes the exceptional corona. Sections~\ref{sec:reference} and~\ref{sec:global} complete the reference-radius comparison and the global identification. The appendices give the geometric coverage and exact arithmetic of the finite exclusions. The supplementary programs and certificate are part of the computer-assisted proof; the main statements do not rely on numerical optimization.

\section{Definitions and radius comparison}\label{sec:definitions}
Let $G=(V,E)$ be the tangency graph of $\PP$. Thus $vw\in E$ precisely when
\[
|x_v-x_w|=r_v+r_w,
\]
where $|\cdot|$ denotes Euclidean norm. The straight-line drawing of $G$ is assumed to be a locally finite triangulation of $\R^2$. In particular, $G$ is connected and every edge belongs to two distinct triangular faces. We write $w\sim v$ for adjacency.

\begin{definition}[Link and link path]\label{def:link}
Let $\mathcal T$ be a triangulation of a surface, possibly with boundary, and let $v$ be a vertex. The \emph{link} $\operatorname{Lk}_{\mathcal T}(v)$ is the graph whose vertices are the neighbors of $v$, with an edge $ab$ exactly when $vab$ is a triangular face of $\mathcal T$. The vertex $v$ itself is not part of its link.

A \emph{link path} is an ordered sequence $(a_1,\ldots,a_k)$ of distinct consecutive vertices in this graph. Its $k-1$ edges correspond to the consecutive triangular faces $va_i a_{i+1}$, $1\le i<k$. These triangles, together with their edges and vertices, form the corresponding \emph{fan} at $v$. Throughout the paper, $k$ counts the vertices of a link path, not its edges; its graph-theoretic length is $k-1$.
\end{definition}

At an interior vertex of a triangulated surface the link is a cycle; at a boundary vertex of a triangulated disk it is a path. Thus in the whole-plane triangulation every link is a finite simple cycle with $d(v)$ vertices. The link records triangular face incidences, rather than all edges between neighbors: an edge $ab$ belongs to the link only if $vab$ is a face.

Scaling all radii and centers by the same positive number preserves $q(\PP)$. We may therefore normalize $\sup_v r_v=1$. A prescribed lower ratio bound $q(\PP)\ge q>0$ then gives $r_v\in[q,1]$ for every $v$. Distinct centers are separated by at least $2q$, and every edge has length at most $2$.

\subsection{Angles and flowers}
For three mutually tangent disks of radii $a,b,c$, let $T(a,b,c)$ be the angle at the center of the $a$-disk in their triangle of centers. The cosine rule, or its half-angle form, gives
\begin{equation}\label{eq:angle}
T(a,b,c)=2\arcsin\sqrt{\frac{bc}{(a+b)(a+c)}}.
\end{equation}
The angle is scale invariant, strictly decreasing in $a$, and strictly increasing in $b$ and $c$.

\begin{definition}
The \emph{flower} at a vertex $v$ is the central disk together with its complete cyclic list of neighbors. If that list is $v_1,\ldots,v_{d(v)}$, its angle sum under a positive radius assignment $\mathbf a=(a_v)_{v\in V}$ is
\begin{equation}\label{eq:angle-sum}
A_v(\mathbf a)=\sum_{i=1}^{d(v)}T(a_v,a_{v_i},a_{v_{i+1}}),
\qquad v_{d(v)+1}=v_1.
\end{equation}
\end{definition}
The actual radii satisfy $A_v(\mathbf r)=2\pi$ at every vertex. Formula~\eqref{eq:angle-sum} also makes sense for an arbitrary positive assignment, whether or not that assignment is realized by an embedded packing. This distinction will be used both for finite exclusions and for the reference assignment.

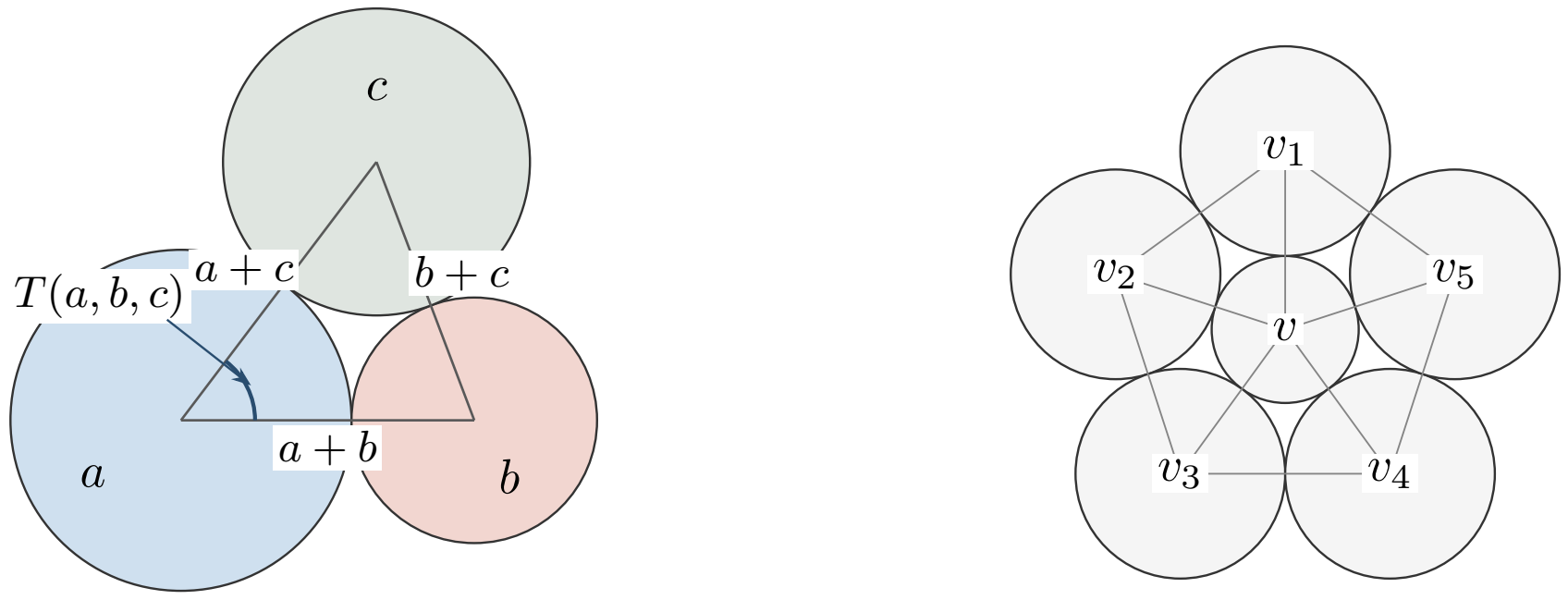
\begin{figure}[htbp]
\centering
\begin{tikzpicture}[font=\small]
\draw[disk,typeL] (1.40000000,1.00000000) circle[radius=1.22000000];
\draw[disk,typeS] (3.49840000,1.00000000) circle[radius=0.87840000];
\draw[disk,typeM] (2.79874419,2.84841519) circle[radius=1.09800000];
\draw[edge] (1.40000000,1.00000000) -- (3.49840000,1.00000000);
\draw[edge] (1.40000000,1.00000000) -- (2.79874419,2.84841519);
\draw[edge] (3.49840000,1.00000000) -- (2.79874419,2.84841519);
\draw[figureaccent,line width=.8pt] (1.93000000,1.00000000) arc[start angle=0,end angle=52.884253,radius=.53];
\node[font=\footnotesize,fill=white,inner sep=1pt] at (0.82000000,1.89000000) {$T(a,b,c)$};
\draw[draw=figureaccent,line width=.45pt,-{Stealth[length=1.4mm]}] (1.30000000,1.72000000) -- (1.90141538,1.24936443);
\node[font=\small] at (0.77000000,0.60000000) {$a$};
\node[font=\small] at (3.74840000,0.60000000) {$b$};
\node[font=\small] at (2.79874419,3.36841519) {$c$};
\node[font=\footnotesize,fill=white,inner sep=1pt] at (2.44920000,0.80000000) {$a+b$};
\node[font=\footnotesize,fill=white,inner sep=1pt] at (1.84937209,2.07420759) {$a+c$};
\node[font=\footnotesize,fill=white,inner sep=1pt] at (3.40857209,2.05420759) {$b+c$};
\node[panel] at (2.45000000,-0.75000000) {(a) A triangle of centers};
\draw[disk,fill=black!4] (9.30000000,1.65000000) circle[radius=0.52597621];
\draw[disk,fill=black!4] (9.30000000,2.92597621) circle[radius=0.75000000];
\draw[disk,fill=black!4] (8.08647451,2.04429833) circle[radius=0.75000000];
\draw[disk,fill=black!4] (8.55000000,0.61771356) circle[radius=0.75000000];
\draw[disk,fill=black!4] (10.05000000,0.61771356) circle[radius=0.75000000];
\draw[disk,fill=black!4] (10.51352549,2.04429833) circle[radius=0.75000000];
\draw[contact] (9.30000000,1.65000000) -- (9.30000000,2.92597621);
\draw[contact] (9.30000000,2.92597621) -- (8.08647451,2.04429833);
\draw[contact] (9.30000000,1.65000000) -- (8.08647451,2.04429833);
\draw[contact] (8.08647451,2.04429833) -- (8.55000000,0.61771356);
\draw[contact] (9.30000000,1.65000000) -- (8.55000000,0.61771356);
\draw[contact] (8.55000000,0.61771356) -- (10.05000000,0.61771356);
\draw[contact] (9.30000000,1.65000000) -- (10.05000000,0.61771356);
\draw[contact] (10.05000000,0.61771356) -- (10.51352549,2.04429833);
\draw[contact] (9.30000000,1.65000000) -- (10.51352549,2.04429833);
\draw[contact] (10.51352549,2.04429833) -- (9.30000000,2.92597621);
\node[font=\small,fill=white,inner sep=1pt] at (9.30000000,1.65000000) {$v$};
\node[font=\small,fill=white,inner sep=1pt] at (9.30000000,2.92597621) {$v_1$};
\node[font=\small,fill=white,inner sep=1pt] at (8.08647451,2.04429833) {$v_2$};
\node[font=\small,fill=white,inner sep=1pt] at (8.55000000,0.61771356) {$v_3$};
\node[font=\small,fill=white,inner sep=1pt] at (10.05000000,0.61771356) {$v_4$};
\node[font=\small,fill=white,inner sep=1pt] at (10.51352549,2.04429833) {$v_5$};
\node[panel] at (9.30000000,-0.75000000) {(b) A complete flower};
\node[annotation] at (9.30000000,-1.20000000) {$\sum_i T(r_v,r_{v_i},r_{v_{i+1}})=2\pi$};
\end{tikzpicture}
\caption{(a) Three mutually tangent disks give a triangle with side lengths $a+b,a+c,b+c$ and angle $T(a,b,c)$ at the first center. (b) The consecutive triangles in a complete flower fill one turn. The right-hand drawing is an illustrative five-flower; it does not specify the degrees of its neighboring vertices.}\label{fig:angles-flower}
\end{figure}

\subsection{The logarithmic comparison formula}
The angle sums behave like a weighted Laplacian in logarithmic radii.

\begin{lemma}\label{lem:jacobian}
Let $\mathbf a$ and $\mathbf b$ be positive radius assignments on a locally finite triangulation without boundary. Set $h_v=\log(a_v/b_v)$. There are positive symmetric weights $w_{vu}=w_{uv}$ on the edges such that
\begin{equation}\label{eq:comparison}
A_v(\mathbf a)-A_v(\mathbf b)
=\sum_{u\sim v}w_{vu}(h_u-h_v).
\end{equation}
If all entries of both assignments lie in one compact subinterval of $(0,\infty)$, the weights have uniform positive lower and finite upper bounds.
\end{lemma}
\begin{proof}
Set $\rho_v(t)=b_v\exp(th_v)$ for $0\le t\le1$.
If $uvz_1$ and $uvz_2$ are the two faces incident to $uv$, define
\begin{equation}\label{eq:edge-weights}
w_{uv}:=\int_0^1
\frac{\sqrt{\rho_u(t)\rho_v(t)}}{\rho_u(t)+\rho_v(t)}
\sum_{j=1}^2
\sqrt{\frac{\rho_{z_j}(t)}
{\rho_u(t)+\rho_v(t)+\rho_{z_j}(t)}}\,dt.
\end{equation}
The verification of \eqref{eq:comparison} is left to the reader;
positivity, symmetry, and the uniform bounds follow from the
formula and compactness.
\end{proof}

For positive symmetric edge weights $w=(w_{uv})$, define
\begin{equation}\label{eq:laplacian}
\Delta_w f(v):=\sum_{u\sim v}w_{vu}(f(u)-f(v)).
\end{equation}
A function is \emph{harmonic} if $\Delta_w f=0$ and \emph{superharmonic} if $\Delta_w f\le0$. The associated random walk moves from $v$ to a neighbor $u$ with probability
\[
p(v,u)=\frac{w_{vu}}{\sum_{z\sim v}w_{vz}}.
\]
All logarithms in this paper are natural. The whole-plane comparison will use recurrence of this walk, not an assumption that $h$ attains a maximum.

\Needspace{10\baselineskip}
\section{Local restrictions}\label{sec:local}
We begin with the published degree restrictions used by the finite exclusions.

\begin{proposition}[Connelly--Zhang]\label{prop:CZ}
Let $\PP$ be a triangulated packing with all radii in $[q,1]$.
\begin{enumerate}
\item If $q\ge0.62$, every vertex has degree $5,6$, or $7$.
\item If $q\ge0.64$, two degree-five vertices cannot be adjacent, and every degree-five vertex has a degree-seven neighbor.
\item If $q\ge0.645$, two triangular faces sharing an edge cannot have all four vertices of degree seven.
\end{enumerate}
\end{proposition}
These are Propositions~7.2--7.5 of~\cite{cz}. We use them as published inputs.

\Needspace{10\baselineskip}
Recall from \eqref{eq:threshold} that $q_0=0.65105$. The following additional restrictions are computer-assisted.

\begin{proposition}[Local exclusions]\label{prop:local}
Let $\PP$ be a triangulated packing with all radii in $[q_0,1]$.
\begin{enumerate}
\item Every degree-five vertex has cyclic neighbor degrees $7,6,7,7,6$. If $ab$ joins its unique consecutive pair of degree-seven neighbors, the vertex in the other face at $ab$ also has degree five.
\item No degree-six vertex has three degree-five neighbors.
\end{enumerate}
\end{proposition}
\begin{proof}
For (i), Proposition~\ref{prop:CZ} reduces the neighbors of a degree-five vertex to degrees six and seven, with at least one of degree seven. The sixteen rooted assignments and their necessary extensions are covered in Appendix~\ref{app:local}. Exact interval tests leave only the stated word and force its opposite degree-five vertex. For (ii), three nonadjacent degree-five neighbors of a degree-six vertex would alternate with degree-seven neighbors. Completing the flowers forced by (i) gives the excluded finite configuration in Appendix~\ref{app:three-small}. The precise configurations, case counts, and validity of the interval tests are given in Appendices~\ref{app:local} and~\ref{app:arithmetic}.
\end{proof}

From now on, unless another bound is stated, assume all radii lie in $[q_0,1]$. The letters $S,M,L$ denote degrees $5,6,7$, not radii. The degree-five word is equivalently
\begin{equation}\label{eq:S-word}
\word{LLMLM}.
\end{equation}

\begin{figure}[htbp]
\centering
\begin{tikzpicture}[font=\small]
\draw[disk,typeL] (-0.88000000,0.00000000) circle[radius=0.88000000];
\draw[disk,typeL] (0.88000000,0.00000000) circle[radius=0.88000000];
\draw[disk,typeS] (0.00000000,1.15610926) circle[radius=0.57292416];
\draw[disk,typeS] (0.00000000,-1.15610926) circle[radius=0.57292416];
\draw[strongedge] (-0.88000000,0.00000000) -- (0.88000000,0.00000000);
\draw[contact] (-0.88000000,0.00000000) -- (0.00000000,1.15610926);
\draw[contact] (0.88000000,0.00000000) -- (0.00000000,1.15610926);
\draw[contact] (-0.88000000,0.00000000) -- (0.00000000,-1.15610926);
\draw[contact] (0.88000000,0.00000000) -- (0.00000000,-1.15610926);
\node[font=\small,fill=white,inner sep=1pt] at (-0.88000000,0.00000000) {$a:L$};
\node[font=\small,fill=white,inner sep=1pt] at (0.88000000,0.00000000) {$b:L$};
\node[font=\small,fill=white,inner sep=1pt] at (0.00000000,1.15610926) {$S$};
\node[font=\small,fill=white,inner sep=1pt] at (0.00000000,-1.15610926) {$S'$};
\node[panel] at (0.00000000,-2.15000000) {(a) The two incident faces};
\draw[twoway] (1.93,1.156109) .. controls (2.70,.7) and (2.70,-.7) .. (1.93,-1.156109);
\node[font=\footnotesize,rotate=90] at (3.00000000,0.00000000) {pairing};
\draw[edge] (7.30000000,0.00000000) -- (6.44771138,1.17307464);
\draw[edge] (6.44771138,1.17307464) -- (8.15228862,1.17307464);
\draw[edge] (7.30000000,0.00000000) -- (8.15228862,1.17307464);
\draw[edge] (8.15228862,1.17307464) -- (8.67903195,-0.44807464);
\draw[edge] (7.30000000,0.00000000) -- (8.67903195,-0.44807464);
\draw[edge] (8.67903195,-0.44807464) -- (7.30000000,-1.45000000);
\draw[edge] (7.30000000,0.00000000) -- (7.30000000,-1.45000000);
\draw[edge] (7.30000000,-1.45000000) -- (5.92096805,-0.44807464);
\draw[edge] (7.30000000,0.00000000) -- (5.92096805,-0.44807464);
\draw[edge] (5.92096805,-0.44807464) -- (6.44771138,1.17307464);
\node[vertex,typeL,minimum size=5.8mm] at (6.44771138,1.17307464) {$\mathsf{L}$};
\node[vertex,typeL,minimum size=5.8mm] at (8.15228862,1.17307464) {$\mathsf{L}$};
\node[vertex,typeM,minimum size=5.8mm] at (8.67903195,-0.44807464) {$\mathsf{M}$};
\node[vertex,typeL,minimum size=5.8mm] at (7.30000000,-1.45000000) {$\mathsf{L}$};
\node[vertex,typeM,minimum size=5.8mm] at (5.92096805,-0.44807464) {$\mathsf{M}$};
\node[vertex,typeS,minimum size=5.8mm] at (7.30000000,0.00000000) {$\mathsf{S}$};
\draw[strongedge] (6.44771138,1.17307464) -- (8.15228862,1.17307464);
\node[panel] at (7.30000000,-2.15000000) {(b) The unique consecutive $LL$ pair};
\end{tikzpicture}
\caption{The pairing of degree-five vertices. At each $S$, the word $LLMLM$ contains a unique consecutive $LL$ pair. Crossing its edge $ab$ gives the paired vertex $S^{\prime}$; crossing back returns $S$. Panel (a) uses the reference radii to display the two triangular faces, with the other neighbors omitted. Panel (b) records only incidences: equal-sized vertex markers do not represent equal disk radii.}\label{fig:pairing}
\end{figure}
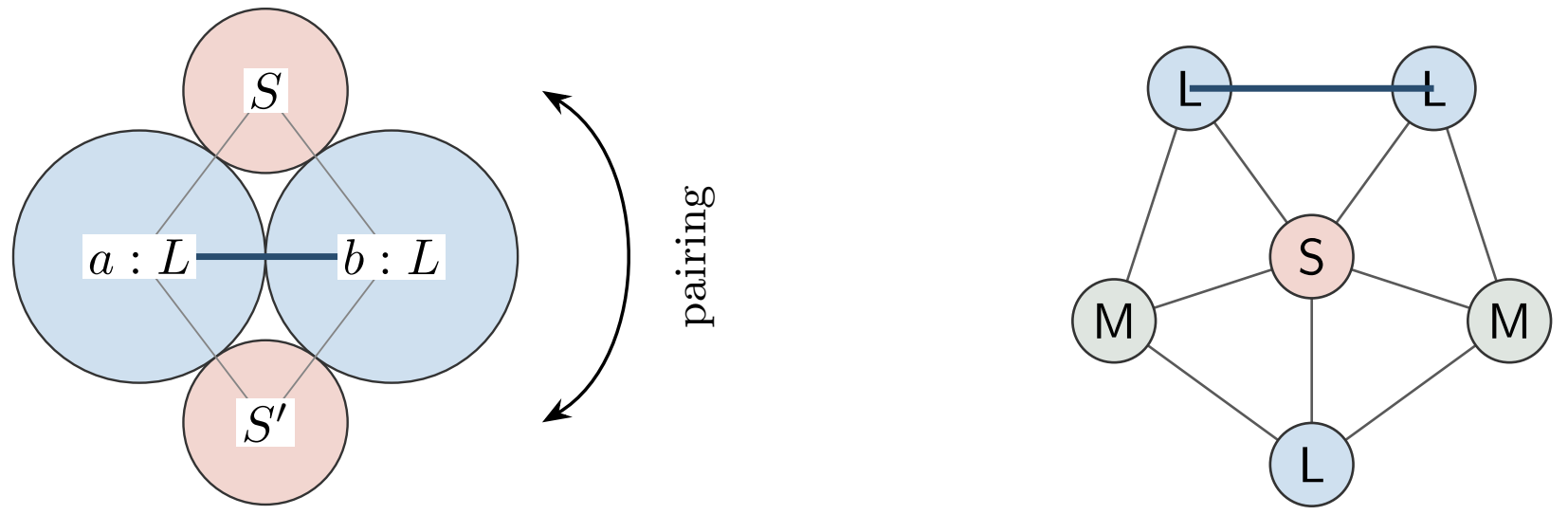

\begin{corollary}[Pairing]\label{cor:pairing}
Crossing the unique consecutive $LL$ edge in the flower of an $S$ vertex defines a fixed-point-free involution on the $S$ vertices.
\end{corollary}
\begin{proof}
Let $v$ be an $S$, and let $ab$ be its consecutive $LL$ edge. The opposite vertex $v'$ is an $S$ by Proposition~\ref{prop:local}(i). In the flower at $v'$, the same edge $ab$ is a consecutive $LL$ pair; it is therefore the unique such pair. Crossing it returns $v$. The opposite vertices in the two faces of an edge are distinct, so $v'\ne v$.
\end{proof}
We call such an edge $ab$ a \emph{paired edge}. The pairing is illustrated in Figure~\ref{fig:pairing}.

\begin{lemma}\label{lem:saturated}
A degree-seven vertex has at most three $S$ neighbors. If it has three, its degree word is one of
\begin{equation}\label{eq:L-two-words}
\word{SLSMSMM},\qquad\word{SMSMSMM}.
\end{equation}
In the first case its only $L$ neighbor is joined to it by a paired edge.
\end{lemma}
\begin{proof}
The $S$ neighbors are independent in the seven-cycle, so there are at most three. When there are three, the numbers of intervening vertices are $1,1,2$.

No $L$ lies in the gap of length two. Indeed, such a neighbor and the central $L$ would form an $LL$ edge incident to an $S$ on one side and a non-$S$ on the other. This contradicts Proposition~\ref{prop:local}(i). The two singleton gaps cannot both contain $L$: their common $S$ would then have two consecutive $LL$ pairs in its own flower, contrary to \eqref{eq:S-word}. Thus either one singleton gap contains $L$ or neither does. These are precisely the two words in \eqref{eq:L-two-words}. In the first case the $LL$ edge has $S$ on both sides and is paired.
\end{proof}

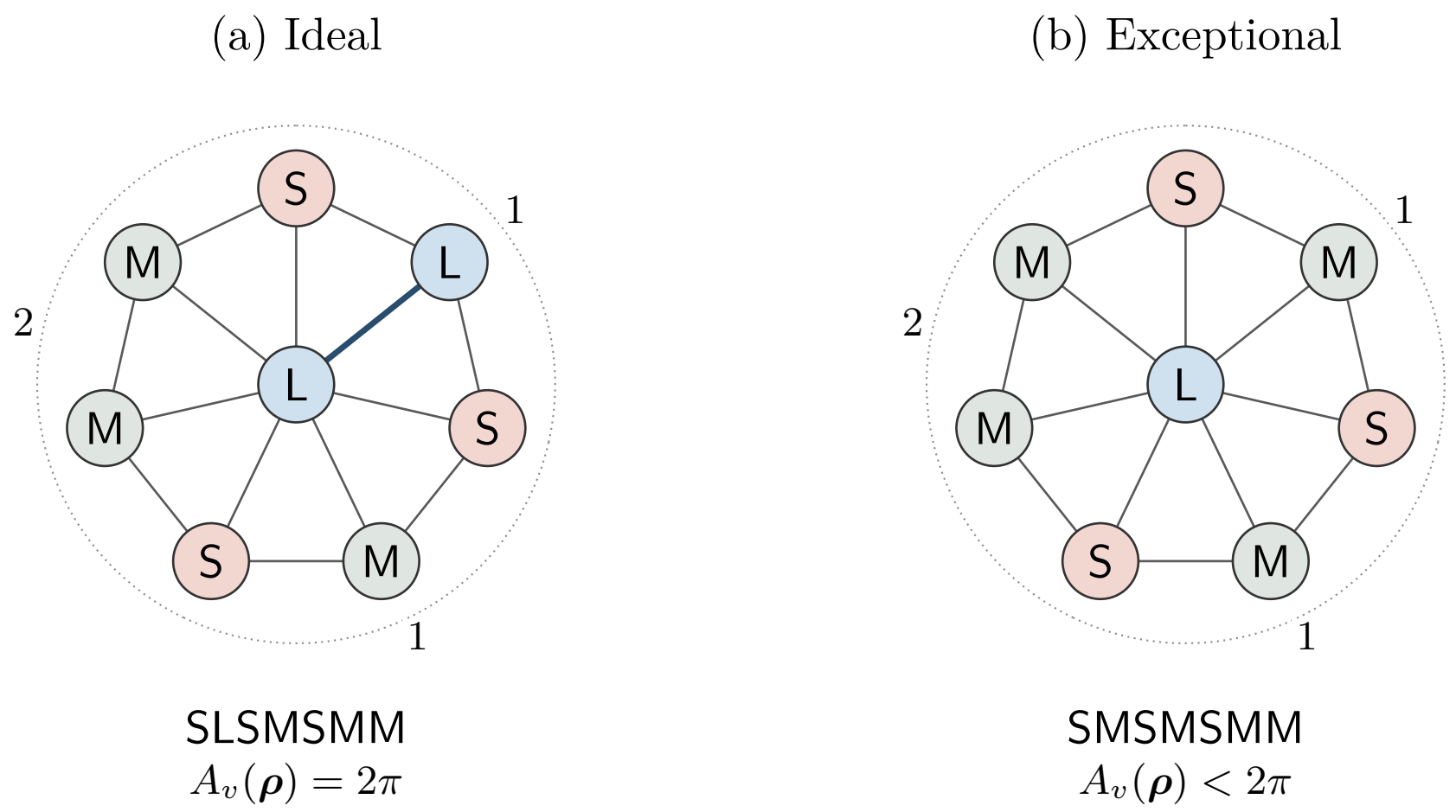
\begin{figure}[htbp]
\centering
\begin{tikzpicture}[font=\small]
\draw[edge] (0.00000000,0.00000000) -- (0.00000000,1.50000000);
\draw[edge] (0.00000000,1.50000000) -- (1.17274722,0.93523470);
\draw[strongedge] (0.00000000,0.00000000) -- (1.17274722,0.93523470);
\draw[edge] (1.17274722,0.93523470) -- (1.46239187,-0.33378140);
\draw[edge] (0.00000000,0.00000000) -- (1.46239187,-0.33378140);
\draw[edge] (1.46239187,-0.33378140) -- (0.65082561,-1.35145330);
\draw[edge] (0.00000000,0.00000000) -- (0.65082561,-1.35145330);
\draw[edge] (0.65082561,-1.35145330) -- (-0.65082561,-1.35145330);
\draw[edge] (0.00000000,0.00000000) -- (-0.65082561,-1.35145330);
\draw[edge] (-0.65082561,-1.35145330) -- (-1.46239187,-0.33378140);
\draw[edge] (0.00000000,0.00000000) -- (-1.46239187,-0.33378140);
\draw[edge] (-1.46239187,-0.33378140) -- (-1.17274722,0.93523470);
\draw[edge] (0.00000000,0.00000000) -- (-1.17274722,0.93523470);
\draw[edge] (-1.17274722,0.93523470) -- (0.00000000,1.50000000);
\node[vertex,typeS,minimum size=5.8mm] at (0.00000000,1.50000000) {$\mathsf{S}$};
\node[vertex,typeL,minimum size=5.8mm] at (1.17274722,0.93523470) {$\mathsf{L}$};
\node[vertex,typeS,minimum size=5.8mm] at (1.46239187,-0.33378140) {$\mathsf{S}$};
\node[vertex,typeM,minimum size=5.8mm] at (0.65082561,-1.35145330) {$\mathsf{M}$};
\node[vertex,typeS,minimum size=5.8mm] at (-0.65082561,-1.35145330) {$\mathsf{S}$};
\node[vertex,typeM,minimum size=5.8mm] at (-1.46239187,-0.33378140) {$\mathsf{M}$};
\node[vertex,typeM,minimum size=5.8mm] at (-1.17274722,0.93523470) {$\mathsf{M}$};
\node[vertex,typeL,minimum size=5.8mm] at (0.00000000,0.00000000) {$\mathsf{L}$};
\draw[draw=black!45,densely dotted] (0.00000000,1.98000000) arc[start angle=90,end angle=-12.857142857,radius=1.98];
\node[font=\footnotesize,fill=white,inner sep=1pt] at (1.67311937,1.33426818) {1};
\draw[draw=black!45,densely dotted] (1.93035727,-0.44059145) arc[start angle=-12.857142857,end angle=-115.714285714,radius=1.98];
\node[font=\footnotesize,fill=white,inner sep=1pt] at (0.92851120,-1.92807338) {1};
\draw[draw=black!45,densely dotted] (-0.85908980,-1.78391836) arc[start angle=-115.714285714,end angle=-270,radius=1.98];
\node[font=\footnotesize,fill=white,inner sep=1pt] at (-2.08634573,0.47619480) {2};
\node[panel] at (0.00000000,-2.62000000) {$\word{SLSMSMM}$};
\node[annotation] at (0.00000000,-3.05000000) {$A_v(\boldsymbol\rho)=2\pi$};
\node[panel] at (0.00000000,2.65000000) {(a) Ideal};
\draw[edge] (6.80000000,0.00000000) -- (6.80000000,1.50000000);
\draw[edge] (6.80000000,1.50000000) -- (7.97274722,0.93523470);
\draw[edge] (6.80000000,0.00000000) -- (7.97274722,0.93523470);
\draw[edge] (7.97274722,0.93523470) -- (8.26239187,-0.33378140);
\draw[edge] (6.80000000,0.00000000) -- (8.26239187,-0.33378140);
\draw[edge] (8.26239187,-0.33378140) -- (7.45082561,-1.35145330);
\draw[edge] (6.80000000,0.00000000) -- (7.45082561,-1.35145330);
\draw[edge] (7.45082561,-1.35145330) -- (6.14917439,-1.35145330);
\draw[edge] (6.80000000,0.00000000) -- (6.14917439,-1.35145330);
\draw[edge] (6.14917439,-1.35145330) -- (5.33760813,-0.33378140);
\draw[edge] (6.80000000,0.00000000) -- (5.33760813,-0.33378140);
\draw[edge] (5.33760813,-0.33378140) -- (5.62725278,0.93523470);
\draw[edge] (6.80000000,0.00000000) -- (5.62725278,0.93523470);
\draw[edge] (5.62725278,0.93523470) -- (6.80000000,1.50000000);
\node[vertex,typeS,minimum size=5.8mm] at (6.80000000,1.50000000) {$\mathsf{S}$};
\node[vertex,typeM,minimum size=5.8mm] at (7.97274722,0.93523470) {$\mathsf{M}$};
\node[vertex,typeS,minimum size=5.8mm] at (8.26239187,-0.33378140) {$\mathsf{S}$};
\node[vertex,typeM,minimum size=5.8mm] at (7.45082561,-1.35145330) {$\mathsf{M}$};
\node[vertex,typeS,minimum size=5.8mm] at (6.14917439,-1.35145330) {$\mathsf{S}$};
\node[vertex,typeM,minimum size=5.8mm] at (5.33760813,-0.33378140) {$\mathsf{M}$};
\node[vertex,typeM,minimum size=5.8mm] at (5.62725278,0.93523470) {$\mathsf{M}$};
\node[vertex,typeL,minimum size=5.8mm] at (6.80000000,0.00000000) {$\mathsf{L}$};
\draw[draw=black!45,densely dotted] (6.80000000,1.98000000) arc[start angle=90,end angle=-12.857142857,radius=1.98];
\node[font=\footnotesize,fill=white,inner sep=1pt] at (8.47311937,1.33426818) {1};
\draw[draw=black!45,densely dotted] (8.73035727,-0.44059145) arc[start angle=-12.857142857,end angle=-115.714285714,radius=1.98];
\node[font=\footnotesize,fill=white,inner sep=1pt] at (7.72851120,-1.92807338) {1};
\draw[draw=black!45,densely dotted] (5.94091020,-1.78391836) arc[start angle=-115.714285714,end angle=-270,radius=1.98];
\node[font=\footnotesize,fill=white,inner sep=1pt] at (4.71365427,0.47619480) {2};
\node[panel] at (6.80000000,-2.62000000) {$\word{SMSMSMM}$};
\node[annotation] at (6.80000000,-3.05000000) {$A_v(\boldsymbol\rho)<2\pi$};
\node[panel] at (6.80000000,2.65000000) {(b) Exceptional};
\end{tikzpicture}
\caption{The two possibilities in Lemma~\ref{lem:saturated}. Three independent $S$ neighbors leave gaps of lengths $1,1,2$, indicated outside the cycles. In the ideal word one singleton gap contains $L$, and the heavy radial edge is paired. In the exceptional word both singleton gaps contain $M$. The ideal word has reference angle sum $2\pi$ by Lemma~\ref{lem:reference-flat}; replacing its singleton $L$ by $M$ strictly lowers that sum. The exceptional word is excluded in Section~\ref{sec:exceptional}. These are combinatorial diagrams, not disk realizations of the exceptional word.}\label{fig:large-coronas}
\end{figure}

Call the first degree-seven word \emph{ideal} and the second \emph{exceptional}. They are shown in Figure~\ref{fig:large-coronas}. Define the reference assignment by
\begin{equation}\label{eq:reference-assignment}
\rho_v=\begin{cases}s,&d(v)=5,\\r,&d(v)=6,\\1,&d(v)=7,\end{cases}
\qquad \boldsymbol\rho=(\rho_v)_{v\in V}.
\end{equation}
This does not assert that actual radii depend only on degree. Its angle sums, shown in the figure, will be analyzed in Section~\ref{sec:reference}.

\section{Excluding the exceptional corona}\label{sec:exceptional}
This section supplies the finite obstruction needed for the whole-plane result, using the cutoff $q_0$ defined in \eqref{eq:threshold}.

\begin{theorem}\label{thm:no-exceptional}
A triangulated disk packing with all radii in $[q_0,1]$ has no exceptional degree-seven vertex.
\end{theorem}

Only Proposition~\ref{prop:local}(i), the degree-five word and its pairing, is used in the geometric reduction below. Proposition~\ref{prop:local}(ii) will instead be used in the reference-angle comparison. We keep this dependence explicit to avoid using any conclusion of the new exclusion as an input to its own proof.

\subsection{The forced neighborhood}
Suppose that $0$ is exceptional, and label its cyclic neighbors $1,\ldots,7$ with degrees
\[
5,6,5,6,5,6,6.
\]
Complete the flowers at these neighbors in order. Denote the vertices opposite $0$ across $12,23,34,45,56,67,71$ by $8,9,10,11,12,13,14$, respectively. The additional labels $15,16,17,18$ complete the flowers at $2,4,6,7$. Figure~\ref{fig:exceptional-initial} shows the resulting source disk.

The flower at $1$ has cyclic list $(0,2,8,14,7)$. Its known neighbor degrees at $0,2,7$ are $7,6,6$. Proposition~\ref{prop:local}(i) therefore forces $8,14$ to have degree seven. Similarly, the flowers at $3$ and $5$ force $9,10,11,12$ to have degree seven. Complete the flowers at
\[
8,9,10,11,12,14
\]
in that order. Across the three edges
\begin{equation}\label{eq:three-paired-edges}
8\,14,\qquad9\,10,\qquad11\,12,
\end{equation}
the pairing rule supplies the degree-five partners $19,24,29$ of $1,3,5$, respectively. Completing these three flowers introduces the final labels $34,35,36$. Figure~\ref{fig:exceptional-complete} shows the full source disk.

\begin{figure}[htbp]
\centering
\begin{tikzpicture}[font=\small]
\tikzset{pn/.style={circle,draw=black!80,fill=white,inner sep=0pt,minimum size=4.8mm,font=\scriptsize,line width=.45pt}}
\coordinate (p8) at (-1.1267656,3.2608127);
\coordinate (p15) at (-2.7108240,2.1339947);
\coordinate (p9) at (-3.4342149,0.3296485);
\coordinate (p10) at (-3.0672669,-1.5793588);
\coordinate (p16) at (-1.7264833,-2.9869308);
\coordinate (p11) at (0.1624466,-3.4461734);
\coordinate (p12) at (1.9998007,-2.8112803);
\coordinate (p17) at (3.2022323,-1.2838256);
\coordinate (p13) at (3.3879778,0.6512347);
\coordinate (p18) at (2.4980642,2.3795326);
\coordinate (p14) at (0.8150329,3.3523456);
\coordinate (p0) at (-0.1275678,-0.0245210);
\coordinate (p1) at (-0.1275678,1.8064192);
\coordinate (p2) at (-1.5610904,1.1921245);
\coordinate (p3) at (-1.8396026,-0.3536071);
\coordinate (p4) at (-1.0078731,-1.6859289);
\coordinate (p5) at (0.5512356,-1.7249821);
\coordinate (p6) at (1.7293717,-0.6570067);
\coordinate (p7) at (1.3625518,1.2513341);
\draw[edge] (p0) -- (p1);
\draw[edge] (p0) -- (p2);
\draw[edge] (p0) -- (p3);
\draw[edge] (p0) -- (p4);
\draw[edge] (p0) -- (p5);
\draw[edge] (p0) -- (p6);
\draw[edge] (p0) -- (p7);
\draw[edge] (p1) -- (p2);
\draw[edge] (p1) -- (p7);
\draw[edge] (p1) -- (p8);
\draw[edge] (p1) -- (p14);
\draw[edge] (p2) -- (p3);
\draw[edge] (p2) -- (p8);
\draw[edge] (p2) -- (p9);
\draw[edge] (p2) -- (p15);
\draw[edge] (p3) -- (p4);
\draw[edge] (p3) -- (p9);
\draw[edge] (p3) -- (p10);
\draw[edge] (p4) -- (p5);
\draw[edge] (p4) -- (p10);
\draw[edge] (p4) -- (p11);
\draw[edge] (p4) -- (p16);
\draw[edge] (p5) -- (p6);
\draw[edge] (p5) -- (p11);
\draw[edge] (p5) -- (p12);
\draw[edge] (p6) -- (p7);
\draw[edge] (p6) -- (p12);
\draw[edge] (p6) -- (p13);
\draw[edge] (p6) -- (p17);
\draw[edge] (p7) -- (p13);
\draw[edge] (p7) -- (p14);
\draw[edge] (p7) -- (p18);
\draw[strongedge] (p8) -- (p14);
\draw[edge] (p8) -- (p15);
\draw[strongedge] (p9) -- (p10);
\draw[edge] (p9) -- (p15);
\draw[edge] (p10) -- (p16);
\draw[strongedge] (p11) -- (p12);
\draw[edge] (p11) -- (p16);
\draw[edge] (p12) -- (p17);
\draw[edge] (p13) -- (p17);
\draw[edge] (p13) -- (p18);
\draw[edge] (p14) -- (p18);
\node[pn,typeL] at (p0) {0};
\node[pn,typeS] at (p1) {1};
\node[pn,typeM] at (p2) {2};
\node[pn,typeS] at (p3) {3};
\node[pn,typeM] at (p4) {4};
\node[pn,typeS] at (p5) {5};
\node[pn,typeM] at (p6) {6};
\node[pn,typeM] at (p7) {7};
\node[pn,typeL] at (p8) {8};
\node[pn,typeL] at (p9) {9};
\node[pn,typeL] at (p10) {10};
\node[pn,typeL] at (p11) {11};
\node[pn,typeL] at (p12) {12};
\node[pn] at (p13) {13};
\node[pn,typeL] at (p14) {14};
\node[pn] at (p15) {15};
\node[pn] at (p16) {16};
\node[pn] at (p17) {17};
\node[pn] at (p18) {18};
\end{tikzpicture}
\caption{The first eight completed flowers forced by an exceptional center $0$. The colors record the known degrees: $S$ (degree five), $M$ (degree six), and $L$ (degree seven). White vertices have no prescribed degree in this drawing. The heavy edges $8\,14$, $9\,10$, and $11\,12$ are the unique consecutive $LL$ edges in the flowers at $1,3,5$. The pairing rule forces a degree-five vertex on their other side. This is a combinatorial source disk, not a disk realization at the threshold $q_0$.}\label{fig:exceptional-initial}
\end{figure}
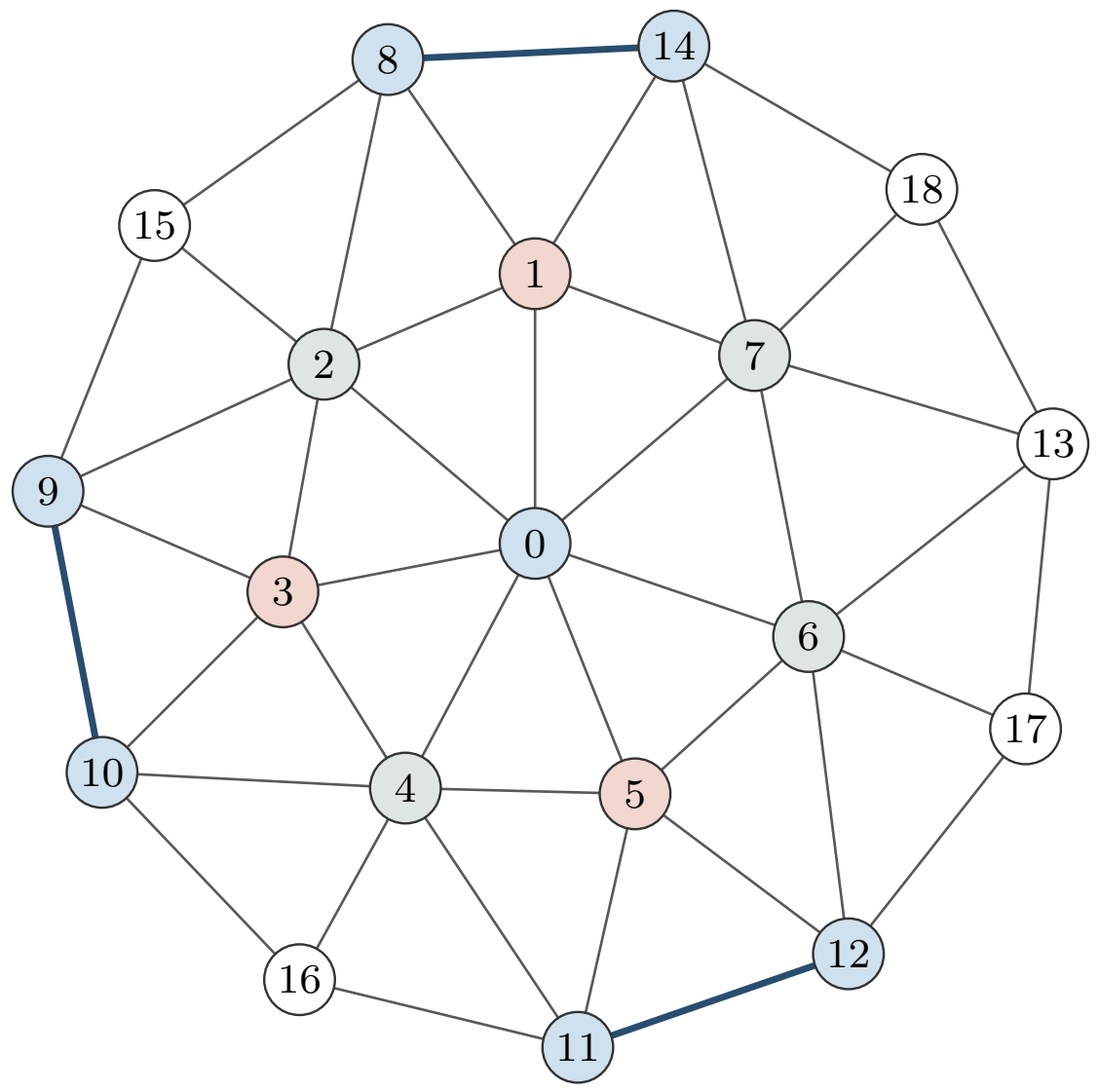
\begin{figure}[p]
\centering
\begin{tikzpicture}[font=\small]
\tikzset{pn/.style={circle,draw=black!80,fill=white,inner sep=0pt,minimum size=4.8mm,font=\scriptsize,line width=.45pt}}
\coordinate (p13) at (5.4230370,1.1803260);
\coordinate (p18) at (4.7928739,2.7983673);
\coordinate (p32) at (3.6935509,4.1424849);
\coordinate (p33) at (2.2326774,5.0811073);
\coordinate (p34) at (0.5532539,5.5223555);
\coordinate (p20) at (-1.1803260,5.4230370);
\coordinate (p21) at (-2.7983673,4.7928739);
\coordinate (p15) at (-4.1424849,3.6935509);
\coordinate (p22) at (-5.0811073,2.2326774);
\coordinate (p23) at (-5.5223555,0.5532539);
\coordinate (p35) at (-5.4230370,-1.1803260);
\coordinate (p26) at (-4.7928739,-2.7983673);
\coordinate (p25) at (-3.6935509,-4.1424849);
\coordinate (p16) at (-2.2326774,-5.0811073);
\coordinate (p27) at (-0.5532539,-5.5223555);
\coordinate (p28) at (1.1803260,-5.4230370);
\coordinate (p36) at (2.7983673,-4.7928739);
\coordinate (p30) at (4.1424849,-3.6935509);
\coordinate (p31) at (5.0811073,-2.2326774);
\coordinate (p17) at (5.5223555,-0.5532539);
\coordinate (p0) at (0.5351482,0.1164752);
\coordinate (p1) at (0.5351482,2.0801272);
\coordinate (p2) at (-1.5606668,1.6426445);
\coordinate (p3) at (-1.6894896,-0.3677180);
\coordinate (p4) at (-0.7367882,-2.1426864);
\coordinate (p5) at (1.3512665,-1.6695484);
\coordinate (p6) at (3.1227405,-0.2801550);
\coordinate (p7) at (2.7238271,1.5526625);
\coordinate (p8) at (-0.9565821,3.6680198);
\coordinate (p9) at (-3.6457409,0.6654116);
\coordinate (p10) at (-3.0394004,-2.1204350);
\coordinate (p11) at (0.6544234,-3.7337848);
\coordinate (p12) at (3.1808086,-2.3075911);
\coordinate (p14) at (1.9340147,3.4208343);
\coordinate (p19) at (0.5166076,4.6230708);
\coordinate (p24) at (-4.4846815,-0.9760925);
\coordinate (p29) at (2.3912820,-3.9901676);
\fill[smallfill!40] (p8) -- (p14) -- (p19) -- cycle;
\fill[smallfill!40] (p8) -- (p19) -- (p20) -- cycle;
\fill[smallfill!40] (p9) -- (p10) -- (p24) -- cycle;
\fill[smallfill!40] (p9) -- (p23) -- (p24) -- cycle;
\fill[smallfill!40] (p10) -- (p24) -- (p26) -- cycle;
\fill[smallfill!40] (p11) -- (p12) -- (p29) -- cycle;
\fill[smallfill!40] (p11) -- (p28) -- (p29) -- cycle;
\fill[smallfill!40] (p12) -- (p29) -- (p30) -- cycle;
\fill[smallfill!40] (p14) -- (p19) -- (p33) -- cycle;
\fill[smallfill!40] (p19) -- (p20) -- (p34) -- cycle;
\fill[smallfill!40] (p19) -- (p33) -- (p34) -- cycle;
\fill[smallfill!40] (p23) -- (p24) -- (p35) -- cycle;
\fill[smallfill!40] (p24) -- (p26) -- (p35) -- cycle;
\fill[smallfill!40] (p28) -- (p29) -- (p36) -- cycle;
\fill[smallfill!40] (p29) -- (p30) -- (p36) -- cycle;
\draw[edge] (p0) -- (p1);
\draw[edge] (p0) -- (p2);
\draw[edge] (p0) -- (p3);
\draw[edge] (p0) -- (p4);
\draw[edge] (p0) -- (p5);
\draw[edge] (p0) -- (p6);
\draw[edge] (p0) -- (p7);
\draw[edge] (p1) -- (p2);
\draw[edge] (p1) -- (p7);
\draw[edge] (p1) -- (p8);
\draw[edge] (p1) -- (p14);
\draw[edge] (p2) -- (p3);
\draw[edge] (p2) -- (p8);
\draw[edge] (p2) -- (p9);
\draw[edge] (p2) -- (p15);
\draw[edge] (p3) -- (p4);
\draw[edge] (p3) -- (p9);
\draw[edge] (p3) -- (p10);
\draw[edge] (p4) -- (p5);
\draw[edge] (p4) -- (p10);
\draw[edge] (p4) -- (p11);
\draw[edge] (p4) -- (p16);
\draw[edge] (p5) -- (p6);
\draw[edge] (p5) -- (p11);
\draw[edge] (p5) -- (p12);
\draw[edge] (p6) -- (p7);
\draw[edge] (p6) -- (p12);
\draw[edge] (p6) -- (p13);
\draw[edge] (p6) -- (p17);
\draw[edge] (p7) -- (p13);
\draw[edge] (p7) -- (p14);
\draw[edge] (p7) -- (p18);
\draw[strongedge] (p8) -- (p14);
\draw[edge] (p8) -- (p15);
\draw[edge] (p8) -- (p19);
\draw[edge] (p8) -- (p20);
\draw[edge] (p8) -- (p21);
\draw[strongedge] (p9) -- (p10);
\draw[edge] (p9) -- (p15);
\draw[edge] (p9) -- (p22);
\draw[edge] (p9) -- (p23);
\draw[edge] (p9) -- (p24);
\draw[edge] (p10) -- (p16);
\draw[edge] (p10) -- (p24);
\draw[edge] (p10) -- (p25);
\draw[edge] (p10) -- (p26);
\draw[strongedge] (p11) -- (p12);
\draw[edge] (p11) -- (p16);
\draw[edge] (p11) -- (p27);
\draw[edge] (p11) -- (p28);
\draw[edge] (p11) -- (p29);
\draw[edge] (p12) -- (p17);
\draw[edge] (p12) -- (p29);
\draw[edge] (p12) -- (p30);
\draw[edge] (p12) -- (p31);
\draw[edge] (p13) -- (p17);
\draw[edge] (p13) -- (p18);
\draw[edge] (p14) -- (p18);
\draw[edge] (p14) -- (p19);
\draw[edge] (p14) -- (p32);
\draw[edge] (p14) -- (p33);
\draw[edge] (p15) -- (p21);
\draw[edge] (p15) -- (p22);
\draw[edge] (p16) -- (p25);
\draw[edge] (p16) -- (p27);
\draw[edge] (p17) -- (p31);
\draw[edge] (p18) -- (p32);
\draw[edge] (p19) -- (p20);
\draw[edge] (p19) -- (p33);
\draw[edge] (p19) -- (p34);
\draw[edge] (p20) -- (p21);
\draw[edge] (p20) -- (p34);
\draw[edge] (p22) -- (p23);
\draw[edge] (p23) -- (p24);
\draw[edge] (p23) -- (p35);
\draw[edge] (p24) -- (p26);
\draw[edge] (p24) -- (p35);
\draw[edge] (p25) -- (p26);
\draw[edge] (p26) -- (p35);
\draw[edge] (p27) -- (p28);
\draw[edge] (p28) -- (p29);
\draw[edge] (p28) -- (p36);
\draw[edge] (p29) -- (p30);
\draw[edge] (p29) -- (p36);
\draw[edge] (p30) -- (p31);
\draw[edge] (p30) -- (p36);
\draw[edge] (p32) -- (p33);
\draw[edge] (p33) -- (p34);
\node[pn,typeL] at (p0) {0};
\node[pn,typeS] at (p1) {1};
\node[pn,typeM] at (p2) {2};
\node[pn,typeS] at (p3) {3};
\node[pn,typeM] at (p4) {4};
\node[pn,typeS] at (p5) {5};
\node[pn,typeM] at (p6) {6};
\node[pn,typeM] at (p7) {7};
\node[pn,typeL] at (p8) {8};
\node[pn,typeL] at (p9) {9};
\node[pn,typeL] at (p10) {10};
\node[pn,typeL] at (p11) {11};
\node[pn,typeL] at (p12) {12};
\node[pn] at (p13) {13};
\node[pn,typeL] at (p14) {14};
\node[pn] at (p15) {15};
\node[pn] at (p16) {16};
\node[pn] at (p17) {17};
\node[pn] at (p18) {18};
\node[pn,typeS] at (p19) {19};
\node[pn] at (p20) {20};
\node[pn] at (p21) {21};
\node[pn] at (p22) {22};
\node[pn] at (p23) {23};
\node[pn,typeS] at (p24) {24};
\node[pn] at (p25) {25};
\node[pn] at (p26) {26};
\node[pn] at (p27) {27};
\node[pn] at (p28) {28};
\node[pn,typeS] at (p29) {29};
\node[pn] at (p30) {30};
\node[pn] at (p31) {31};
\node[pn] at (p32) {32};
\node[pn] at (p33) {33};
\node[pn] at (p34) {34};
\node[pn] at (p35) {35};
\node[pn] at (p36) {36};
\end{tikzpicture}
\caption{The unfolded $37$-label disk used to exclude the exceptional corona. The seventeen completed centers are the colored vertices. The partners of $1,3,5$ are $19,24,29$, across the three heavy $LL$ edges. The shaded fans complete these partners and introduce labels $34,35,36$. White vertices have free boundary radii. Every row of Table~\ref{tab:flowers} is visible in this diagram. Different source labels may have the same image in an actual packing; injectivity of the entire patch is not assumed.}\label{fig:exceptional-complete}
\end{figure}
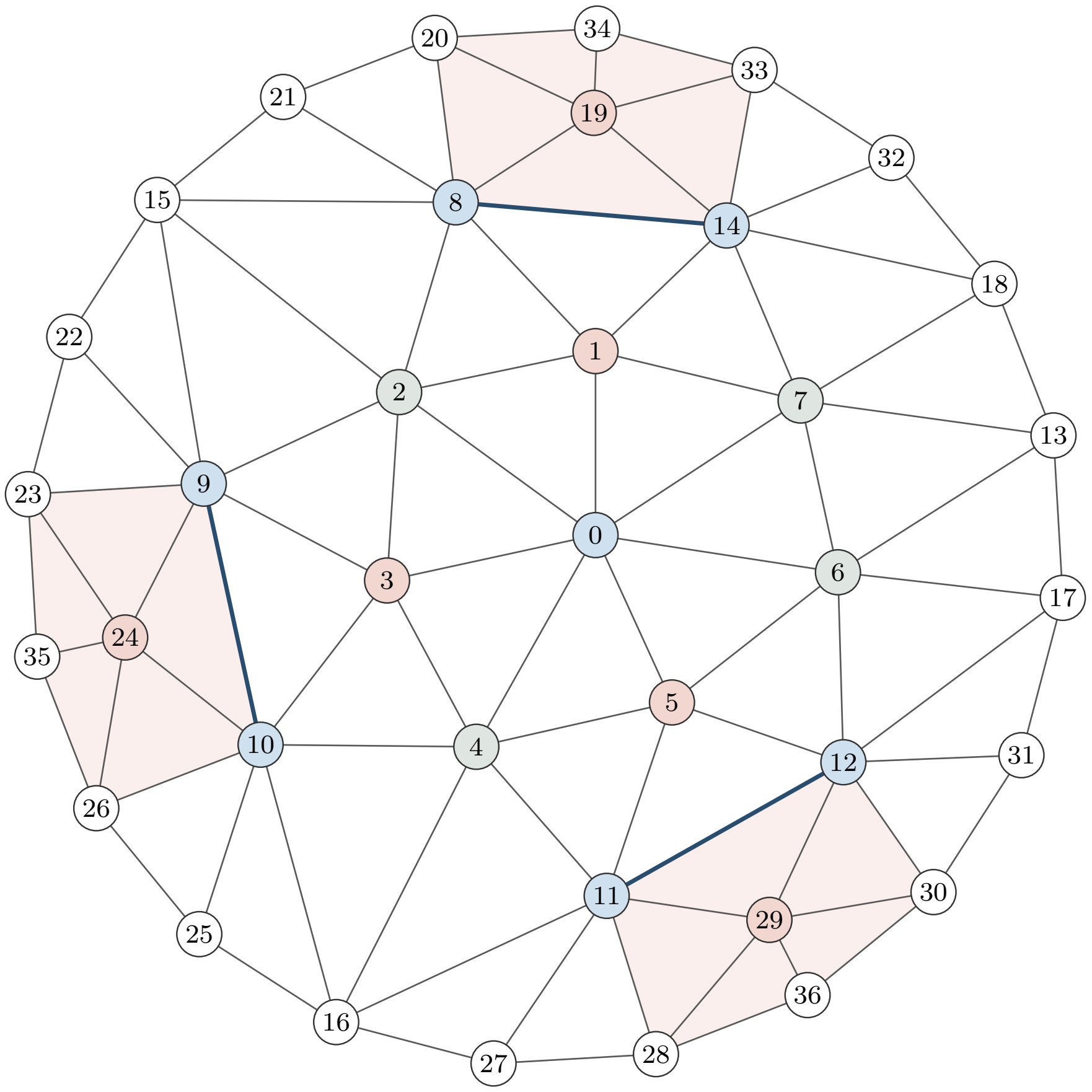

\begin{table}[htbp]
\centering
\caption{The seventeen completed flowers in the exceptional-corona system. Every row is a cyclic list; reversal preserves its angle equation. The labels range from $0$ to $36$, and different source labels need not represent different disks.}\label{tab:flowers}
\begin{tabular}{@{}rcl@{}}
\toprule
Center $i$ & Degree & Cyclic neighbor list $J_i$\\
\midrule
0&7&$1,2,3,4,5,6,7$\\
1&5&$0,2,8,14,7$\\
2&6&$0,1,8,15,9,3$\\
3&5&$0,2,9,10,4$\\
4&6&$0,3,10,16,11,5$\\
5&5&$0,4,11,12,6$\\
6&6&$0,5,12,17,13,7$\\
7&6&$0,1,14,18,13,6$\\
8&7&$1,2,15,21,20,19,14$\\
9&7&$2,3,10,24,23,22,15$\\
10&7&$3,4,16,25,26,24,9$\\
11&7&$4,5,12,29,28,27,16$\\
12&7&$5,6,17,31,30,29,11$\\
14&7&$1,7,18,32,33,19,8$\\
19&5&$8,14,33,34,20$\\
24&5&$9,10,26,35,23$\\
29&5&$11,12,30,36,28$\\
\bottomrule
\end{tabular}
\end{table}

\subsection{Unfolding and possible identifications}
Here the \emph{source} is the abstract finite triangulated disk just constructed, and the \emph{target} is the triangulation of the putative packing. Links and link paths are as in Definition~\ref{def:link}: a boundary vertex of the source has a path as its link, whereas its image in the target has a cycle as its link. The source disk is not assumed to occur injectively in the packing. The following argument permits the geometric reduction even when distant labels have the same image.

\begin{lemma}[Link extension]\label{lem:link}
Let an oriented triangulated disk be mapped face by face into an oriented triangulation of a surface without boundary, with each face mapped to a face and the two source faces at an interior edge mapped to opposite sides of its image. Suppose a boundary vertex has a link path with $k$ vertices (and $k-1$ edges), and its image has degree $d>k$. Its flower can be completed to degree $d$ by introducing $d-k$ fresh source labels and $d-k+1$ triangular faces, and the map extends to the completed fan.
\end{lemma}
\begin{proof}
The target link is a simple $d$-cycle. Consecutive source faces give consecutive edges of this cycle. The opposite-side condition prevents backtracking, so the image of the source path follows one fixed cyclic direction. Since $k<d$, its $k$ image vertices are distinct. The complementary arc between the two endpoint images has $d-k$ interior vertices. Assign these, in order, to the fresh source labels and attach the corresponding $d-k+1$ triangular faces. The completed flower traverses the target link exactly once. Images of the new labels may coincide with images of source vertices elsewhere in the disk, but this does not obstruct the face-by-face map.
\end{proof}

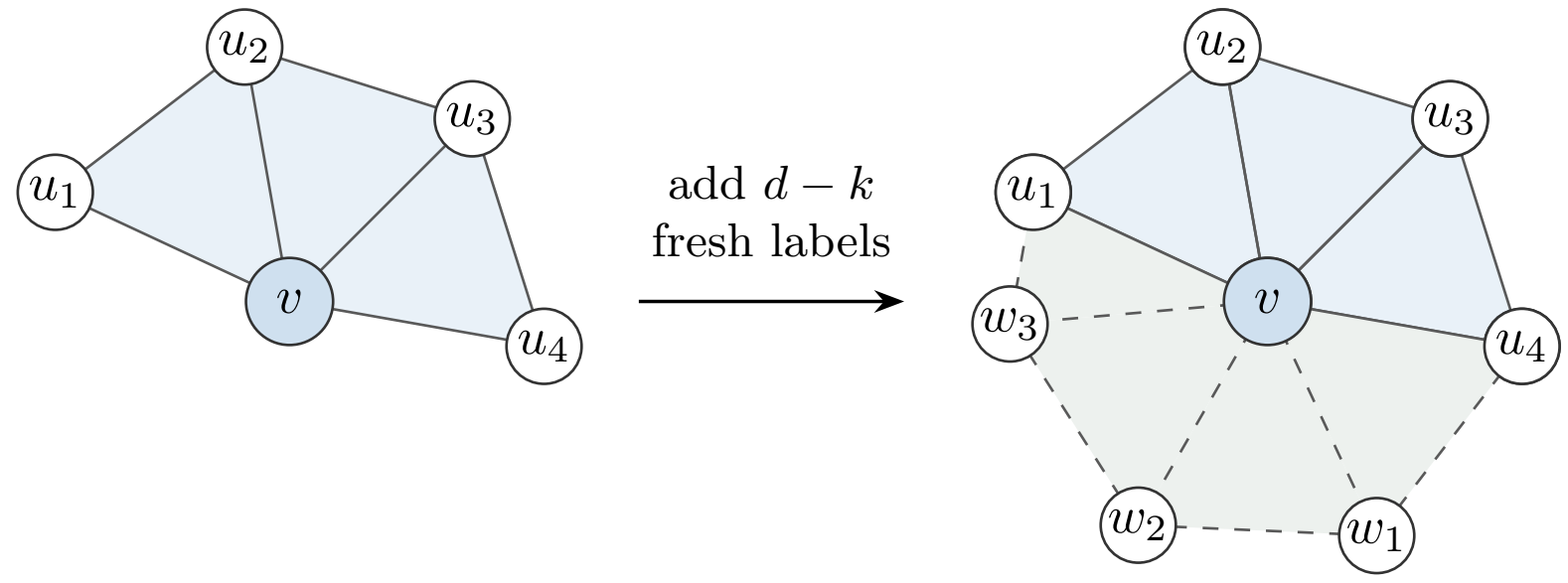
\begin{figure}[htbp]
\centering
\begin{tikzpicture}[font=\small,scale=.93]
\foreach \shift/\title in {0/{(a) A boundary link path},7/{(b) The completed flower}} {
\begin{scope}[xshift=\shift cm]
\coordinate (v) at (0,0);
\foreach \i/\ang in {1/155,2/100,3/45,4/-10,5/-65,6/-120,7/-175} {
\coordinate (u\i) at (\ang:1.85);
}
\foreach \i/\j in {1/2,2/3,3/4}{
\fill[largefill!45] (v)--(u\i)--(u\j)--cycle;
\draw[edge] (u\i)--(u\j);
}
\foreach \i in {1,2,3,4}{\draw[edge] (v)--(u\i);\node[vertex,minimum size=5mm] at (u\i) {$u_{\i}$};}
\node[vertex,typeL] at (v) {$v$};
\node[panel] at (0,-2.65) {\title};
\end{scope}
}
\begin{scope}[xshift=7cm]
\foreach \i/\j in {4/5,5/6,6/7,7/1}{
\fill[mediumfill!55] (v)--(u\i)--(u\j)--cycle;
\draw[edge,dashed] (u\i)--(u\j);
}
\foreach \i/\j in {5/1,6/2,7/3}{\draw[edge,dashed] (v)--(u\i);\node[vertex,fill=white,minimum size=5mm] at (u\i) {$w_{\j}$};}
\foreach \i in {1,2,3,4}{\draw[edge] (v)--(u\i);\node[vertex,minimum size=5mm] at (u\i) {$u_{\i}$};}
\node[vertex,typeL] at (v) {$v$};
\end{scope}
\draw[arrow] (2.5,0)--(4.4,0);
\node[annotation] at (3.45,.65) {add $d-k$\\fresh labels};
\end{tikzpicture}
\caption{Link extension, shown for $k=4$ and $d=7$; see Definition~\ref{def:link}. The $u_i$ are the $k$ existing link vertices, joined by $k-1$ link edges. The $w_i$ are the $d-k$ fresh source labels. Completing the link cycle adds $d-k+1$ triangular faces at $v$. The new labels are distinct in the source disk, although their images may coincide with vertices elsewhere in the target.}\label{fig:link-extension}
\end{figure}

Let $\II$ be the set of seventeen centers in Table~\ref{tab:flowers}. For $i\in\II$, write $J_i=(j_1,\ldots,j_{d_i})$, where $d_i=|J_i|$, set $j_{d_i+1}=j_1$, and define
\begin{equation}\label{eq:finite-angle}
\calA_i(\mathbf a)=\sum_{k=1}^{d_i}T(a_i,a_{j_k},a_{j_{k+1}}),
\qquad \mathbf a=(a_0,\ldots,a_{36}).
\end{equation}
Consider the finite system
\begin{equation}\label{eq:finite-system}
\mathbf a\in[q_0,1]^{37},
\qquad \calA_i(\mathbf a)=2\pi\quad(i\in\II).
\end{equation}
No equations are imposed at the other twenty vertices. There are no nonedge separation constraints and no assumed equalities between radius variables.

\begin{lemma}[Geometric reduction]\label{lem:geometric-reduction}
An exceptional vertex in a packing with all radii in $[q_0,1]$ supplies a solution of \eqref{eq:finite-system}.
\end{lemma}
\begin{proof}
Start with the actual flower at $0$ and use Lemma~\ref{lem:link} in the completion order above. For the first seven completions, the numbers $k$ of vertices in the existing link paths and the target degrees $d$ are
\[
\begin{array}{c|rrrrrrr}
\text{center}&1&2&3&4&5&6&7\\\hline
k&3&4&4&4&4&4&5\\
d&5&6&5&6&5&6&6.
\end{array}
\]
For the remaining completions, they are
\[
\begin{array}{c|rrrrrrrrr}
\text{center}&8&9&10&11&12&14&19&24&29\\\hline
k&4&4&5&4&5&5&4&4&4\\
d&7&7&7&7&7&7&5&5&5.
\end{array}
\]
Every completion has $k<d$. The degrees in the first table are given by the exceptional word; those in the second are forced by Proposition~\ref{prop:local}(i), as explained above. Thus every row of Table~\ref{tab:flowers} maps to one complete actual flower. The source disk has $37$ vertices, $88$ edges, $52$ triangular faces, and $20$ boundary vertices; the displayed lists also determine these counts directly.

Assign to each source label the radius of its image disk. These radii belong to $[q_0,1]$ and satisfy every equation in \eqref{eq:finite-system}. If labels have a common image, the corresponding radii are equal. Omitting those equalities only enlarges the feasible set.
\end{proof}

\subsection{Exclusion of the finite angle system}
The decisive step is a global exclusion of the box in \eqref{eq:finite-system}, not a local optimization of one realization.

\begin{proposition}[Computer-assisted finite exclusion]\label{prop:finite-exclusion}
The system \eqref{eq:finite-system} has no solution.
\end{proposition}
\begin{proof}
The exact interval rules are established in Appendix~\ref{app:arithmetic}. For a box of radius intervals, monotonicity gives a lower and upper bound for each complete angle sum. A box is rejected only when one of these bounds excludes $2\pi$. Radius intervals can be contracted using the same monotonicity, and a subdivision uses both closed children, so no point is lost at a splitting endpoint.

The supplied finite certificate starts with $[65105/100000,1]^{37}$. It has $3894$ binary splits, $3895$ rejected leaves, and $834176$ verified radius-bound contractions. In all, $7789$ boxes are processed and $838071$ strict angle-sum inequalities are verified, with no unresolved or pending box. A separate integer-only checker reconstructs every box and validates every angular enclosure by the cosine inequalities in Lemma~\ref{lem:cosine-enclosure}. It hard-codes Table~\ref{tab:flowers} and the initial interval, rather than trusting these data to the proof trace. The certificate-covering lemma, Lemma~\ref{lem:certificate-soundness}, then proves infeasibility.
\end{proof}

\begin{proof}[Proof of Theorem~\ref{thm:no-exceptional}]
By Lemma~\ref{lem:geometric-reduction}, an exceptional vertex would supply a solution of \eqref{eq:finite-system}. Proposition~\ref{prop:finite-exclusion} excludes such a solution.
\end{proof}

\begin{remark}
The exceptional-corona exclusion at $q_0$ uses the complete extended system. The seven triangles at the exceptional center alone do not give a contradiction. Completing the three paired degree-five flowers is essential to the obstruction certified here. The drawings in this section specify incidences, not a claimed packing at the forbidden radius bound.
\end{remark}

\section{The reference radius assignment}\label{sec:reference}
We now compare arbitrary high-ratio packings with the reference packing. Recall the degree-based reference assignment $\boldsymbol\rho$ from \eqref{eq:reference-assignment}. We do not yet assume that it is proportional to the actual radii.

The reference small and medium flowers in packing~53 have radius words $11r1r$ and $1r1s1s$, respectively~\cite{fhs,fhs-short}. Consequently,
\begin{align}
T(s,1,1)+4T(s,1,r)&=2\pi,\label{eq:small-reference}\\
2T(r,1,r)+4T(r,1,s)&=2\pi.\label{eq:medium-reference}
\end{align}
Equation~\eqref{eq:small-reference}, together with the isolated value of $s$, characterizes $r$ in $(0.83420,0.83440)$. The published construction supplies these exact radii and attainment of the final bound. The fixed-radius uniqueness recalled in the introduction also gives a short route to the final identification; Section~\ref{subsec:identification} proves it directly from the local words. Optimality over arbitrary radius sets is established by the present comparison argument.

\begin{lemma}\label{lem:reference-flat}
For the reference assignment, all three words in \eqref{eq:ideal-words} have angle sum $2\pi$.
\end{lemma}
\begin{proof}
The first two words give \eqref{eq:small-reference} and \eqref{eq:medium-reference}. Denote their sums by $A_S,A_M$, and let $A_L$ be the sum for $\word{SLSMSMM}$ centered at radius $1$. Thus
\[
A_L=2T(1,s,1)+4T(1,s,r)+T(1,r,r).
\]
The sum $A_S+A_M+A_L$ groups into the three angles of one triangle with radii $s,1,1$, four triangles with radii $s,r,1$, and one triangle with radii $1,r,r$. It is therefore $6\pi$. Since $A_S=A_M=2\pi$, also $A_L=2\pi$.
\end{proof}

\begin{figure}[htbp]
\centering
\begin{tikzpicture}[font=\small]
\draw[disk,typeS] (0.00000000,0.00000000) circle[radius=0.42969312];
\draw[disk,typeL] (0.64050555,0.88158025) circle[radius=0.66000000];
\draw[disk,typeL] (-0.67916639,0.85215264) circle[radius=0.66000000];
\draw[disk,typeM] (-0.92170333,-0.33394595) circle[radius=0.55064199];
\draw[disk,typeL] (0.02429324,-1.08942230) circle[radius=0.66000000];
\draw[disk,typeM] (0.93567319,-0.29252795) circle[radius=0.55064199];
\draw[contact] (0.00000000,0.00000000) -- (0.64050555,0.88158025);
\draw[contact] (0.64050555,0.88158025) -- (-0.67916639,0.85215264);
\draw[contact] (0.00000000,0.00000000) -- (-0.67916639,0.85215264);
\draw[contact] (-0.67916639,0.85215264) -- (-0.92170333,-0.33394595);
\draw[contact] (0.00000000,0.00000000) -- (-0.92170333,-0.33394595);
\draw[contact] (-0.92170333,-0.33394595) -- (0.02429324,-1.08942230);
\draw[contact] (0.00000000,0.00000000) -- (0.02429324,-1.08942230);
\draw[contact] (0.02429324,-1.08942230) -- (0.93567319,-0.29252795);
\draw[contact] (0.00000000,0.00000000) -- (0.93567319,-0.29252795);
\draw[contact] (0.93567319,-0.29252795) -- (0.64050555,0.88158025);
\node[font=\small,fill=white,fill opacity=.83,text opacity=1,inner sep=1pt] at (0.00000000,0.00000000) {$\mathsf{S}$};
\node[font=\small,fill=white,fill opacity=.83,text opacity=1,inner sep=1pt] at (0.64050555,0.88158025) {$\mathsf{L}$};
\node[font=\small,fill=white,fill opacity=.83,text opacity=1,inner sep=1pt] at (-0.67916639,0.85215264) {$\mathsf{L}$};
\node[font=\small,fill=white,fill opacity=.83,text opacity=1,inner sep=1pt] at (-0.92170333,-0.33394595) {$\mathsf{M}$};
\node[font=\small,fill=white,fill opacity=.83,text opacity=1,inner sep=1pt] at (0.02429324,-1.08942230) {$\mathsf{L}$};
\node[font=\small,fill=white,fill opacity=.83,text opacity=1,inner sep=1pt] at (0.93567319,-0.29252795) {$\mathsf{M}$};
\node[panel] at (0.00000000,-2.37000000) {$d(v)=5,\quad\rho_v=s$};
\node[annotation] at (0.00000000,-2.82000000) {$\word{LLMLM}$};
\draw[disk,typeM] (4.90000000,0.00000000) circle[radius=0.55064199];
\draw[disk,typeL] (4.90000000,1.21064199) circle[radius=0.66000000];
\draw[disk,typeM] (3.91922306,0.50090217) circle[radius=0.55064199];
\draw[disk,typeL] (3.91922306,-0.70973982) circle[radius=0.66000000];
\draw[disk,typeS] (4.97552683,-0.97742142) circle[radius=0.42969312];
\draw[disk,typeL] (5.97816855,-0.55064199) circle[radius=0.66000000];
\draw[disk,typeS] (5.73611571,0.51182755) circle[radius=0.42969312];
\draw[contact] (4.90000000,0.00000000) -- (4.90000000,1.21064199);
\draw[contact] (4.90000000,1.21064199) -- (3.91922306,0.50090217);
\draw[contact] (4.90000000,0.00000000) -- (3.91922306,0.50090217);
\draw[contact] (3.91922306,0.50090217) -- (3.91922306,-0.70973982);
\draw[contact] (4.90000000,0.00000000) -- (3.91922306,-0.70973982);
\draw[contact] (3.91922306,-0.70973982) -- (4.97552683,-0.97742142);
\draw[contact] (4.90000000,0.00000000) -- (4.97552683,-0.97742142);
\draw[contact] (4.97552683,-0.97742142) -- (5.97816855,-0.55064199);
\draw[contact] (4.90000000,0.00000000) -- (5.97816855,-0.55064199);
\draw[contact] (5.97816855,-0.55064199) -- (5.73611571,0.51182755);
\draw[contact] (4.90000000,0.00000000) -- (5.73611571,0.51182755);
\draw[contact] (5.73611571,0.51182755) -- (4.90000000,1.21064199);
\node[font=\small,fill=white,fill opacity=.83,text opacity=1,inner sep=1pt] at (4.90000000,0.00000000) {$\mathsf{M}$};
\node[font=\small,fill=white,fill opacity=.83,text opacity=1,inner sep=1pt] at (4.90000000,1.21064199) {$\mathsf{L}$};
\node[font=\small,fill=white,fill opacity=.83,text opacity=1,inner sep=1pt] at (3.91922306,0.50090217) {$\mathsf{M}$};
\node[font=\small,fill=white,fill opacity=.83,text opacity=1,inner sep=1pt] at (3.91922306,-0.70973982) {$\mathsf{L}$};
\node[font=\small,fill=white,fill opacity=.83,text opacity=1,inner sep=1pt] at (4.97552683,-0.97742142) {$\mathsf{S}$};
\node[font=\small,fill=white,fill opacity=.83,text opacity=1,inner sep=1pt] at (5.97816855,-0.55064199) {$\mathsf{L}$};
\node[font=\small,fill=white,fill opacity=.83,text opacity=1,inner sep=1pt] at (5.73611571,0.51182755) {$\mathsf{S}$};
\node[panel] at (4.90000000,-2.37000000) {$d(v)=6,\quad\rho_v=r$};
\node[annotation] at (4.90000000,-2.82000000) {$\word{LMLSLS}$};
\draw[disk,typeL] (9.80000000,0.00000000) circle[radius=0.66000000];
\draw[disk,typeS] (9.80000000,1.08969312) circle[radius=0.42969312];
\draw[disk,typeL] (8.74966010,0.79949114) circle[radius=0.66000000];
\draw[disk,typeS] (8.74966010,-0.29020198) circle[radius=0.42969312];
\draw[disk,typeM] (9.29904525,-1.10213355) circle[radius=0.55064199];
\draw[disk,typeS] (10.27201136,-0.98215904) circle[radius=0.42969312];
\draw[disk,typeM] (10.97354664,-0.29739252) circle[radius=0.55064199];
\draw[disk,typeM] (10.72891913,0.77637830) circle[radius=0.55064199];
\draw[contact] (9.80000000,0.00000000) -- (9.80000000,1.08969312);
\draw[contact] (9.80000000,1.08969312) -- (8.74966010,0.79949114);
\draw[contact] (9.80000000,0.00000000) -- (8.74966010,0.79949114);
\draw[contact] (8.74966010,0.79949114) -- (8.74966010,-0.29020198);
\draw[contact] (9.80000000,0.00000000) -- (8.74966010,-0.29020198);
\draw[contact] (8.74966010,-0.29020198) -- (9.29904525,-1.10213355);
\draw[contact] (9.80000000,0.00000000) -- (9.29904525,-1.10213355);
\draw[contact] (9.29904525,-1.10213355) -- (10.27201136,-0.98215904);
\draw[contact] (9.80000000,0.00000000) -- (10.27201136,-0.98215904);
\draw[contact] (10.27201136,-0.98215904) -- (10.97354664,-0.29739252);
\draw[contact] (9.80000000,0.00000000) -- (10.97354664,-0.29739252);
\draw[contact] (10.97354664,-0.29739252) -- (10.72891913,0.77637830);
\draw[contact] (9.80000000,0.00000000) -- (10.72891913,0.77637830);
\draw[contact] (10.72891913,0.77637830) -- (9.80000000,1.08969312);
\node[font=\small,fill=white,fill opacity=.83,text opacity=1,inner sep=1pt] at (9.80000000,0.00000000) {$\mathsf{L}$};
\node[font=\small,fill=white,fill opacity=.83,text opacity=1,inner sep=1pt] at (9.80000000,1.08969312) {$\mathsf{S}$};
\node[font=\small,fill=white,fill opacity=.83,text opacity=1,inner sep=1pt] at (8.74966010,0.79949114) {$\mathsf{L}$};
\node[font=\small,fill=white,fill opacity=.83,text opacity=1,inner sep=1pt] at (8.74966010,-0.29020198) {$\mathsf{S}$};
\node[font=\small,fill=white,fill opacity=.83,text opacity=1,inner sep=1pt] at (9.29904525,-1.10213355) {$\mathsf{M}$};
\node[font=\small,fill=white,fill opacity=.83,text opacity=1,inner sep=1pt] at (10.27201136,-0.98215904) {$\mathsf{S}$};
\node[font=\small,fill=white,fill opacity=.83,text opacity=1,inner sep=1pt] at (10.97354664,-0.29739252) {$\mathsf{M}$};
\node[font=\small,fill=white,fill opacity=.83,text opacity=1,inner sep=1pt] at (10.72891913,0.77637830) {$\mathsf{M}$};
\node[panel] at (9.80000000,-2.37000000) {$d(v)=7,\quad\rho_v=1$};
\node[annotation] at (9.80000000,-2.82000000) {$\word{SLSMSMM}$};
\end{tikzpicture}
\caption{The three reference flowers, drawn at a common scale with radii $s,r,1$. In each panel $v$ is the central vertex; the neighbor labels are their degree types in the reference packing. Each flower closes with angle sum $2\pi$. Lemma~\ref{lem:reference-flat} uses only these local words and their angle identities, not a classification of the full triangulation.}\label{fig:reference-flowers}
\end{figure}

A degree-six vertex with word $\word{MMMMMM}$ will be called \emph{passive}; one with word $\word{LMLSLS}$ will be called \emph{active}. These names refer only to degree words.

\begin{proposition}[Reference-angle comparison]\label{prop:reference-excess}
In a triangulated packing with all radii in $[q_0,1]$,
\begin{equation}\label{eq:reference-excess}
A_v(\boldsymbol\rho)\ge2\pi
\qquad\text{for every }v.
\end{equation}
Equality holds exactly for the three words in \eqref{eq:ideal-words} and for the passive degree-six word.
\end{proposition}
\begin{proof}
At an $S$ vertex the word is \eqref{eq:S-word}, so Lemma~\ref{lem:reference-flat} applies.

At an $M$ vertex there are at most two $S$ neighbors, by Proposition~\ref{prop:local}(ii). Each such $S$ is flanked by $L$ in the $M$ flower: the two $M$ slots in the $S$ word $\word{LLMLM}$ are each flanked by $L$, and these incidences describe the same triangular faces.

If there is no $S$, every neighbor has type $M$ or $L$, and hence reference radius $r$ or $1$. The central reference radius is $r$, so monotonicity gives angle sum at least $6T(r,r,r)=2\pi$, with equality only for the passive word. This is a statement about the auxiliary assignment $\boldsymbol\rho$; the actual radii have not yet been determined. If there is exactly one $S$, lowering all unspecified neighbors to $r$ gives the radius word $s1rrr1$. If there are two opposite $S$ slots, their flanks force $s11s11$. If two $S$ slots are at cyclic distance two, their three distinct flanks are $L$, leaving one slot. A remaining $M$ gives the active word and equality; a remaining $L$ strictly increases that angle sum.

At an $L$ vertex with at most two $S$ neighbors, lower the remaining radii to $r$, and if necessary choose extra nonadjacent slots to lower to $s$. The resulting word has two nonadjacent $s$ slots and five $r$ slots. Its angle sum is
\[
4T(1,s,r)+3T(1,r,r),
\]
independently of the separation of the two slots. If there are three $S$ neighbors, Lemma~\ref{lem:saturated} and Theorem~\ref{thm:no-exceptional} leave only the ideal $L$ word, which gives equality.

The three remaining strict comparisons are recorded in Table~\ref{tab:reference-bounds}. They follow from the rational bounds
\begin{equation}\label{eq:reference-bounds}
0.65105018<s<0.65105019,
\qquad 0.83420<r<0.83440
\end{equation}
and the exact angle inequalities of Appendix~\ref{app:arithmetic}. Each lower bound in the table is greater than $2\pi<6.283186$. This proves the assertion and the equality cases.
\end{proof}

\begin{table}[htbp]
\centering
\caption{Strict reference-angle comparisons. In each row the central radius is fixed by its degree, and the listed neighbors give a lower angle sum. The endpoints are rational certified bounds, not floating-point assertions.}\label{tab:reference-bounds}
\small
\begin{tabular}{@{}lll@{}}
\toprule
Case & Lower angle sum & Certified lower bound\\
\midrule
$M$, one $S$ & $2T(r,s,1)+2T(r,1,r)+2T(r,r,r)$ & $6.333966$\\
$M$, two opposite $S$ & $4T(r,s,1)+2T(r,1,1)$ & $6.391892$\\
$L$, at most two $S$ & $4T(1,s,r)+3T(1,r,r)$ & $6.331199$\\
\bottomrule
\end{tabular}
\end{table}

\section{Comparison and identification on the whole plane}\label{sec:global}
The comparison first determines the radii and local words. We then identify the entire packing. Both steps are needed: the exceptional-corona exclusion alone does not force every degree-seven vertex to have three $S$ neighbors.

\subsection{Recurrence and the comparison principle}
A bounded ratio between two radius assignments on an infinite packing need not attain its infimum. Recurrence supplies the minimum principle for bounded superharmonic functions that is needed below.

Dubejko's comparison theorem~\cite[Theorem~3.1]{dubejko} applies to positive labels on a recurrent triangulation, including auxiliary labels that need not be realized by an embedded packing. The argument below gives a self-contained proof of the bounded-radius case used here. For the general relation between recurrence and circle-packing type, see He and Schramm~\cite{he-schramm}.

\begin{lemma}\label{lem:recurrence}
Let $G$ be the tangency graph of a triangulated packing with radii in $[a,1]$, where $a>0$, and suppose its degrees are bounded. Give its edges positive symmetric weights bounded above by one finite constant. The associated random walk is recurrent, and every bounded superharmonic function is constant.
\end{lemma}
\begin{proof}
Fix a disk center as origin. For all sufficiently large integers $n$, let $\mathcal C_n$ be the set of edges joining a center in the closed Euclidean disk of radius $3n$ to a center outside it. These are finite cutsets separating the origin from infinity. They are pairwise edge-disjoint, because an edge has length at most $2$.

The endpoints of $\mathcal C_n$ lie in an annulus of bounded width about the circle of radius $3n$. Disks of radius $a$ about the centers have disjoint interiors, so an area estimate bounds the number of these endpoints by a constant times $n$. Bounded degrees and the upper weight bound then give
\[
w(\mathcal C_n):=\sum_{\{u,v\}\in\mathcal C_n}w_{uv}\le Cn
\]
for some fixed $C<\infty$. Thus $\sum_n1/w(\mathcal C_n)=\infty$, and the Nash--Williams criterion~\cite{lyonsperes} gives recurrence.

Let $f$ be bounded and superharmonic. Start the conductance walk $(X_n)_{n\ge0}$ at a vertex $v$, and let $\tau_u$ be its first hitting time of a fixed vertex $u$. The process $f(X_n)$ is a bounded supermartingale. Let $\mathbb E_v$ denote expectation for this walk. Stopping at $\min\{\tau_u,N\}$ gives $f(v)\ge\mathbb E_v[f(X_{\min\{\tau_u,N\}})]$. Recurrence implies $\tau_u<\infty$ almost surely. Bounded convergence as $N\to\infty$ yields $f(v)\ge f(u)$. Interchanging $u,v$ proves equality.
\end{proof}

The next lemma also follows from the six-circle theorem of B\'ar\'any, F\"uredi, and Pach~\cite[Theorem~1]{bfp}. They proved that a packing in which every disk has at least six neighbors is the regular equal-disk packing if its radii have positive infimum. In the degree range of Proposition~\ref{prop:CZ}, a nonuniform packing must therefore have a degree-five disk. We give a proof using the comparison formula.

\begin{lemma}\label{lem:degree-five}
A nonuniform triangulated packing with all radii in $[0.62,1]$ contains a degree-five vertex.
\end{lemma}
\begin{proof}
Compare the actual radii with the constant assignment $1$ by Lemma~\ref{lem:jacobian}. For $u_v=\log r_v$, the actual angle sum $2\pi$ and the reference sum $d(v)\pi/3$ give
\begin{equation}\label{eq:degree-charge}
\Delta_w u(v)=(6-d(v))\frac\pi3.
\end{equation}
By Proposition~\ref{prop:CZ}, the degrees are $5,6,7$. In the absence of degree five, $u$ is bounded and superharmonic. Lemma~\ref{lem:recurrence} makes it constant, contrary to nonuniformity.
\end{proof}

\subsection{Identification from the ideal words}\label{subsec:identification}
The following lemmas determine the global triangulation from the three degree words in
\eqref{eq:ideal-words}.  The argument is combinatorial.  In particular, the
regular polygons used below supply an auxiliary metric; they are not asserted
to be polygons of disk centers in the original packing.

Replacing forced disk clusters by polygonal cells is also useful in Kennedy's
binary classifications~\cite[Section~3]{kennedy}. Here the three ideal words
determine both the resulting cell decomposition and its restoration.

We regard a plane triangulation as a cell decomposition of the plane, with
its vertices, edges, triangular faces, and their incidences. An isomorphism
of plane cell decompositions consists of incidence-preserving bijections
between their vertices, edges, and faces; it is allowed to reverse orientation.
Two vertices of an octagon are called opposite if they are four edges apart
along its boundary. The letters $S,M,L$ continue to denote vertices of
degrees $5,6,7$, respectively.

\begin{lemma}\label{lem:quadrilateral-octagon-reduction}
Let $\mathcal T$ be a locally finite triangulation of the plane with the
degree words in \eqref{eq:ideal-words}.  Deleting all $S$ vertices and all
$LL$ and $MM$ edges gives a locally finite cell decomposition $\mathcal D$
of the plane with the following properties:
\begin{enumerate}
\item Every face is a quadrilateral or an octagon with a simple boundary,
and its boundary vertices alternate between $L$ and $M$.
\item Every vertex has degree three and belongs to one quadrilateral and
two distinct octagons.
\item The original triangulation is recovered by putting the $MM$
diagonal into each quadrilateral and choosing one of the two diagonals
joining opposite $L$ vertices in each octagon, then inserting an $S$
vertex joined to the five boundary vertices of each resulting pentagon.
\end{enumerate}
\end{lemma}
\begin{proof}
The words show that every triangular face has type $SLL$, $SLM$, or $LMM$.
Every $L$ has exactly one $L$ neighbor.  If $uv$ is an $LL$ edge, the two
faces incident to it have distinct third vertices $a,b$ of type $S$.
Conversely, the word at an $S$ vertex has exactly one consecutive $LL$
pair, so this construction partitions the $S$ vertices into pairs.

Deleting an $S$ vertex leaves a pentagonal disk whose boundary word is
$LLMLM$.  We first show that the pentagons associated with paired vertices
$a,b$ intersect precisely in their common $LL$ edge $uv$.
Let $z_a$ and $z_b$ be the third $L$ neighbors of $a$ and $b$, respectively.
In the link of $z_a$, the entry $a$ is flanked by two $M$ entries, since
$z_a$ is the $L$ entry between the two $M$ entries in the word at $a$.
The word $SLSMSMM$ at an $L$ vertex has exactly one $S$ entry flanked by
two $M$ entries.  Thus $z_a\ne z_b$.

Suppose that $a$ and $b$ have a common $M$ neighbor $m$.  They are then
the two $S$ entries in the word $LMLSLS$ at $m$, and those entries have
a common intervening $L$ neighbor $z$.  The sets of $L$ neighbors of
$a$ and $b$ intersect precisely in $\{u,v\}$, so $z$ is $u$ or $v$.
Write $z'$ for the other endpoint of $uv$.  The link of $z$ contains
both paths $(a,z',b)$ and $(a,m,b)$.  Their four vertices are distinct,
and their union is a cycle of length four.  This is impossible in the
simple link cycle of length seven at $z$.  Hence $a$ and $b$ have no
common $M$ neighbor.

The two pentagonal boundaries therefore have only $u,v$ in common.
Their interiors consist of distinct triangular faces of $\mathcal T$,
so their union is a disk.  Removing $uv$ from this union gives an
octagon with eight distinct boundary vertices.  Its boundary alternates
between $L$ and $M$, and $u,v$ are opposite $L$ vertices.

Every $M$ has exactly one $M$ neighbor.  The two faces incident to an
$MM$ edge have type $LMM$ and have distinct third vertices.  Removing
that edge merges them into a quadrilateral with simple alternating
boundary $LMLM$.
Every original triangular face belongs to exactly one of the octagonal
or quadrilateral disks just described.  Their interiors are disjoint,
so they form a cell decomposition of the plane.  Local finiteness is
inherited from $\mathcal T$, since each new face is a union of at most
ten original triangular faces.

An $L$ vertex loses its three edges to $S$ vertices and its unique
$LL$ edge, leaving degree three.  Its two $S$ neighbors next to its
$L$ neighbor form one pair, and its third $S$ neighbor belongs to
another pair.  It consequently belongs to two distinct octagons.
It is incident to one $LMM$ triangle and therefore to one quadrilateral.
An $M$ vertex loses its two edges to $S$ vertices and its unique
$MM$ edge, again leaving degree three.  Its two $S$ neighbors cannot
be paired, by the preceding common-neighbor argument, so it also belongs
to two distinct octagons and one quadrilateral.

Finally, each deleted $MM$ edge is the diagonal joining the two $M$
vertices of its quadrilateral.  Each deleted $LL$ edge joins opposite
$L$ vertices of its octagon and divides it into the two original
pentagons. Inserting one interior $S$ vertex in each pentagon and joining
it to the five boundary vertices recovers the two $S$ flowers. These operations recover every face
and edge of $\mathcal T$.
\end{proof}

\begin{figure}[htbp]
\centering
\begin{tikzpicture}[font=\small]
\path[fill=smallfill!35] (2.380000000,0.000000000) -- (2.049030663,0.799030663) -- (1.250000000,1.130000000) -- (0.450969337,0.799030663) -- (0.120000000,0.000000000) -- cycle;
\path[fill=smallfill!18] (0.120000000,0.000000000) -- (0.450969337,-0.799030663) -- (1.250000000,-1.130000000) -- (2.049030663,-0.799030663) -- (2.380000000,0.000000000) -- cycle;
\path[edge,fill=none] (2.380000000,0.000000000) -- (2.049030663,0.799030663) -- (1.250000000,1.130000000) -- (0.450969337,0.799030663) -- (0.120000000,0.000000000) -- (0.450969337,-0.799030663) -- (1.250000000,-1.130000000) -- (2.049030663,-0.799030663) -- cycle;
\path[edge,fill=none] (6.230000000,0.000000000) -- (5.899030663,0.799030663) -- (5.100000000,1.130000000) -- (4.300969337,0.799030663) -- (3.970000000,0.000000000) -- (4.300969337,-0.799030663) -- (5.100000000,-1.130000000) -- (5.899030663,-0.799030663) -- cycle;
\draw[contact,dashed] (1.250000000,0.440000000) -- (2.380000000,0.000000000);
\draw[contact,dashed] (1.250000000,0.440000000) -- (2.049030663,0.799030663);
\draw[contact,dashed] (1.250000000,0.440000000) -- (1.250000000,1.130000000);
\draw[contact,dashed] (1.250000000,0.440000000) -- (0.450969337,0.799030663);
\draw[contact,dashed] (1.250000000,0.440000000) -- (0.120000000,0.000000000);
\node[vertex,typeS,minimum size=4.3mm,font=\footnotesize] at (1.250000000,0.440000000) {$ S $};
\draw[contact,dashed] (1.250000000,-0.440000000) -- (0.120000000,0.000000000);
\draw[contact,dashed] (1.250000000,-0.440000000) -- (0.450969337,-0.799030663);
\draw[contact,dashed] (1.250000000,-0.440000000) -- (1.250000000,-1.130000000);
\draw[contact,dashed] (1.250000000,-0.440000000) -- (2.049030663,-0.799030663);
\draw[contact,dashed] (1.250000000,-0.440000000) -- (2.380000000,0.000000000);
\node[vertex,typeS,minimum size=4.3mm,font=\footnotesize] at (1.250000000,-0.440000000) {$ S $};
\draw[strongedge,dashed] (2.380000000,0.000000000) -- (0.120000000,0.000000000);
\node[vertex,typeL,minimum size=4.4mm,font=\footnotesize] at (2.380000000,0.000000000) {$ L $};
\node[vertex,typeM,minimum size=4.4mm,font=\footnotesize] at (2.049030663,0.799030663) {$ M $};
\node[vertex,typeL,minimum size=4.4mm,font=\footnotesize] at (1.250000000,1.130000000) {$ L $};
\node[vertex,typeM,minimum size=4.4mm,font=\footnotesize] at (0.450969337,0.799030663) {$ M $};
\node[vertex,typeL,minimum size=4.4mm,font=\footnotesize] at (0.120000000,0.000000000) {$ L $};
\node[vertex,typeM,minimum size=4.4mm,font=\footnotesize] at (0.450969337,-0.799030663) {$ M $};
\node[vertex,typeL,minimum size=4.4mm,font=\footnotesize] at (1.250000000,-1.130000000) {$ L $};
\node[vertex,typeM,minimum size=4.4mm,font=\footnotesize] at (2.049030663,-0.799030663) {$ M $};
\node[vertex,typeL,minimum size=4.4mm,font=\footnotesize] at (6.230000000,0.000000000) {$ L $};
\node[vertex,typeM,minimum size=4.4mm,font=\footnotesize] at (5.899030663,0.799030663) {$ M $};
\node[vertex,typeL,minimum size=4.4mm,font=\footnotesize] at (5.100000000,1.130000000) {$ L $};
\node[vertex,typeM,minimum size=4.4mm,font=\footnotesize] at (4.300969337,0.799030663) {$ M $};
\node[vertex,typeL,minimum size=4.4mm,font=\footnotesize] at (3.970000000,0.000000000) {$ L $};
\node[vertex,typeM,minimum size=4.4mm,font=\footnotesize] at (4.300969337,-0.799030663) {$ M $};
\node[vertex,typeL,minimum size=4.4mm,font=\footnotesize] at (5.100000000,-1.130000000) {$ L $};
\node[vertex,typeM,minimum size=4.4mm,font=\footnotesize] at (5.899030663,-0.799030663) {$ M $};
\draw[arrow] (2.70,0) -- (3.65,0);
\node[annotation] at (3.175,1.72) {delete the $S$'s\\and the $LL$ edge};
\node[panel] at (3.175,-1.70) {(a) Paired five-flowers};
\path[fill=mediumfill!35] (8.250000000,0.960000000) -- (7.290000000,0.000000000) -- (8.250000000,-0.960000000) -- (9.210000000,0.000000000) -- cycle;
\path[edge,fill=none] (8.250000000,0.960000000) -- (7.290000000,0.000000000) -- (8.250000000,-0.960000000) -- (9.210000000,0.000000000) -- cycle;
\path[edge,fill=none] (12.100000000,0.960000000) -- (11.140000000,0.000000000) -- (12.100000000,-0.960000000) -- (13.060000000,0.000000000) -- cycle;
\draw[strongedge,dashed] (7.290000000,0.000000000) -- (9.210000000,0.000000000);
\node[vertex,typeL,minimum size=4.4mm,font=\footnotesize] at (8.250000000,0.960000000) {$ L $};
\node[vertex,typeM,minimum size=4.4mm,font=\footnotesize] at (7.290000000,0.000000000) {$ M $};
\node[vertex,typeL,minimum size=4.4mm,font=\footnotesize] at (8.250000000,-0.960000000) {$ L $};
\node[vertex,typeM,minimum size=4.4mm,font=\footnotesize] at (9.210000000,0.000000000) {$ M $};
\node[vertex,typeL,minimum size=4.4mm,font=\footnotesize] at (12.100000000,0.960000000) {$ L $};
\node[vertex,typeM,minimum size=4.4mm,font=\footnotesize] at (11.140000000,0.000000000) {$ M $};
\node[vertex,typeL,minimum size=4.4mm,font=\footnotesize] at (12.100000000,-0.960000000) {$ L $};
\node[vertex,typeM,minimum size=4.4mm,font=\footnotesize] at (13.060000000,0.000000000) {$ M $};
\draw[arrow] (9.62,0) -- (10.72,0);
\node[annotation] at (10.17,1.72) {delete the\\$MM$ edge};
\node[panel] at (10.17,-1.70) {(b) An $MM$ edge};
\end{tikzpicture}
\caption{The two coarsening operations. In (a), deleting the paired $S$ vertices leaves two pentagons; deleting their common $LL$ edge joins them into one alternating octagon. In (b), deleting an $MM$ edge joins its two incident triangles into an alternating quadrilateral. Solid boundary edges remain and dashed interior edges are removed. The drawings show incidences, not disk-center geometry.}\label{fig:global-uniqueness}
\end{figure}
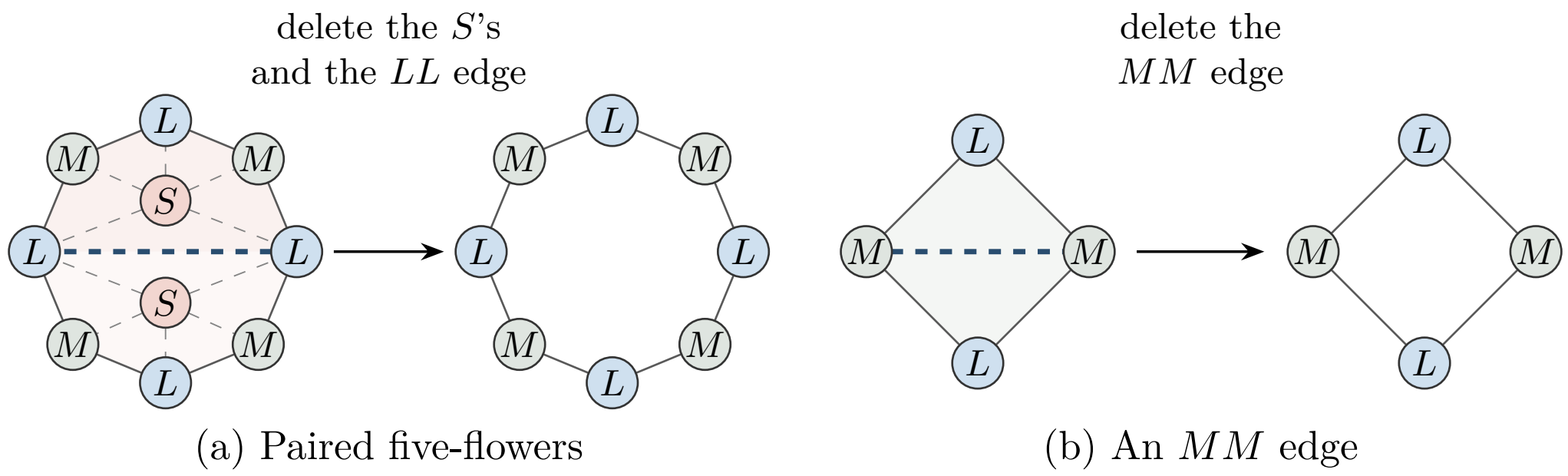

Figure~\ref{fig:global-uniqueness} illustrates the two merging operations.
The notation $4.8.8$ records the resulting incidence of one quadrilateral
and two octagons at every vertex.

\begin{lemma}\label{lem:square-octagon-tiling}
Let $\mathcal D$ be a locally finite cell decomposition of the plane
whose faces have simple boundaries of lengths four and eight.  Suppose
that every vertex has degree three and is incident to one quadrilateral
and two distinct octagons.  Then $\mathcal D$ is isomorphic to the
standard tiling by regular squares and regular octagons in which one
square and two octagons meet at every vertex.
\end{lemma}
\begin{proof}
Give each quadrilateral the metric of a regular square of side length
one, give each octagon the metric of a regular octagon of side length
one, and identify corresponding sides by isometries. Measure path lengths
in these polygons and define distance as the infimum of path lengths.
This gives an intrinsic length metric on the underlying plane. The metric is
Euclidean at face interiors and edge interiors.  At every vertex,
the sum of the incident angles is
\[
\frac{\pi}{2}+\frac{3\pi}{4}+\frac{3\pi}{4}=2\pi,
\]
so the metric is Euclidean there as well.

This metric is complete.  Indeed, the two polygonal shapes and the
vertex configuration are fixed.  Small Euclidean neighborhoods in
face interiors, across edge interiors, and around vertices therefore
give a number $\varepsilon>0$, independent of the point, such that
every point has a metric neighborhood isometric to a Euclidean disk
of radius $\varepsilon$.  A metric Cauchy sequence eventually lies
in the concentric disk of radius $\varepsilon/2$ in one such
neighborhood and hence converges.

The resulting surface is simply connected, since its underlying
topological space is the plane.  Continuing its Euclidean charts
along paths gives a developing map to $\R^2$; simple connectedness
makes this map single valued.  It is a local isometry.  Completeness
allows every straight segment in $\R^2$ to be lifted for its full
length, so the developing map is a covering map onto $\R^2$.
Since $\R^2$ is simply connected, it is an isometry.  Thus the faces
of $\mathcal D$ form an edge-to-edge tiling of the Euclidean plane
by the specified regular polygons.

For completeness, this regular tiling is unique.  Along the boundary
of an octagon, the faces across successive edges alternate between
squares and octagons: at each endpoint, exactly one of the other two
incident faces is a square.  Every edge of a square has an octagon
on its other side.  Once one octagon and the type of the face across
one of its edges are fixed, these rules determine the type and position
of every face by propagation across edges. The \emph{face adjacency graph}
has one vertex for each face and an edge for each pair of faces sharing an
edge. For a locally finite tiling of the plane this graph is connected,
so the propagation determines the whole tiling.
The required tiling exists: start with a square grid of side length
$1+\sqrt2$ and cut off each corner of every grid square at distance
$1/\sqrt2$ along its two sides.  The remaining cells are regular
octagons of side length one, and the gaps at the grid vertices are
regular squares of side length one.  This is the standard tiling
in the statement.
\end{proof}

\begin{proposition}\label{prop:unique-ideal-triangulation}
Every locally finite triangulation of the plane with the degree words
in \eqref{eq:ideal-words} is isomorphic to the tangency triangulation
of packing~53.
\end{proposition}
\begin{proof}
Apply Lemmas~\ref{lem:quadrilateral-octagon-reduction} and
\ref{lem:square-octagon-tiling}.  It remains to show that there is,
up to isomorphism, only one way to recover the original triangulation
from the standard square--octagon tiling.

Every edge of the reduced graph joins an $L$ vertex to an $M$ vertex.
Thus its vertex types give a bipartition. A connected bipartite graph
has exactly two proper colorings by the labels $L,M$, exchanged by
interchanging the labels.
A quarter-turn about the center of a square is a symmetry of the
standard tiling and takes a square vertex to an adjacent vertex;
it therefore interchanges the two parts.  Consequently we may identify
the reduced graphs of any two candidate triangulations with the same
standard tiling and with the same $L,M$ labels.

The $MM$ diagonal in each square is uniquely determined. An octagon
has four $L$ vertices in alternating positions, and there are exactly
two diagonals joining opposite $L$ vertices.  We show that a choice
in one octagon determines all the others; see Figure~\ref{fig:global-propagation}.

An edge separating two octagons has exactly one $L$ endpoint, say $v$.
Those two octagons are precisely the octagons incident to $v$.
In the original triangulation, $v$ has exactly one $L$ neighbor.
It follows that exactly one of the two octagons must select its
diagonal incident to $v$.  In each octagon exactly one of its two
possible diagonals is incident to $v$.  The choice in either octagon
therefore uniquely determines the choice in the other.

The graph whose vertices are the octagons and whose edges join
octagons sharing an edge is connected: in the grid construction in
Lemma~\ref{lem:square-octagon-tiling}, it is the adjacency graph of
the original grid squares.  Hence there are at most two possible
global diagonal assignments, distinguished by the choice in a fixed
octagon.  A quarter-turn about the center of that octagon preserves
the standard tiling, preserves its $L,M$ bipartition, and interchanges
the two choices.  Indeed, it moves each octagon vertex by two positions
along the boundary.  After this symmetry is applied, the propagation
rule shows that any two valid global assignments agree everywhere.

The diagonal and the inserted $S$ vertices determine all the original
triangles by Lemma~\ref{lem:quadrilateral-octagon-reduction}.  There
is therefore at most one plane triangulation with the three given
degree words, up to isomorphism.  Packing~53 has these words and
provides an example, proving the assertion.
\end{proof}

\begin{figure}[htbp]
\centering
\begin{tikzpicture}[font=\small]
\path[fill=mediumfill!50] (2.300182900,1.685000000) -- (2.915365799,2.300182900) -- (2.300182900,2.915365799) -- (1.685000000,2.300182900) -- cycle;
\draw[edge] (2.300182900,1.685000000) -- (1.685000000,2.300182900);
\draw[edge] (0.815000000,2.300182900) -- (1.685000000,2.300182900);
\draw[edge] (0.815000000,2.300182900) -- (0.199817100,1.685000000);
\draw[edge] (0.199817100,1.685000000) -- (0.199817100,0.815000000);
\draw[edge] (0.199817100,0.815000000) -- (0.815000000,0.199817100);
\draw[edge] (0.815000000,0.199817100) -- (1.685000000,0.199817100);
\draw[edge] (1.685000000,0.199817100) -- (2.300182900,0.815000000);
\draw[edge] (2.300182900,1.685000000) -- (2.300182900,0.815000000);
\draw[edge] (4.400548699,1.685000000) -- (3.785365799,2.300182900);
\draw[edge] (2.915365799,2.300182900) -- (3.785365799,2.300182900);
\draw[edge] (2.915365799,2.300182900) -- (2.300182900,1.685000000);
\draw[edge] (2.300182900,0.815000000) -- (2.915365799,0.199817100);
\draw[edge] (2.915365799,0.199817100) -- (3.785365799,0.199817100);
\draw[edge] (3.785365799,0.199817100) -- (4.400548699,0.815000000);
\draw[edge] (4.400548699,0.815000000) -- (4.400548699,1.685000000);
\draw[edge] (2.300182900,3.785365799) -- (1.685000000,4.400548699);
\draw[edge] (1.685000000,4.400548699) -- (0.815000000,4.400548699);
\draw[edge] (0.815000000,4.400548699) -- (0.199817100,3.785365799);
\draw[edge] (0.199817100,3.785365799) -- (0.199817100,2.915365799);
\draw[edge] (0.199817100,2.915365799) -- (0.815000000,2.300182900);
\draw[edge] (1.685000000,2.300182900) -- (2.300182900,2.915365799);
\draw[edge] (2.300182900,3.785365799) -- (2.300182900,2.915365799);
\draw[edge] (4.400548699,3.785365799) -- (3.785365799,4.400548699);
\draw[edge] (3.785365799,4.400548699) -- (2.915365799,4.400548699);
\draw[edge] (2.915365799,4.400548699) -- (2.300182900,3.785365799);
\draw[edge] (2.300182900,2.915365799) -- (2.915365799,2.300182900);
\draw[edge] (3.785365799,2.300182900) -- (4.400548699,2.915365799);
\draw[edge] (4.400548699,2.915365799) -- (4.400548699,3.785365799);
\node[vertex,typeL,minimum size=4.4mm,font=\footnotesize] at (0.199817100,0.815000000) {$ L $};
\node[vertex,typeM,minimum size=4.4mm,font=\footnotesize] at (0.199817100,1.685000000) {$ M $};
\node[vertex,typeM,minimum size=4.4mm,font=\footnotesize] at (0.199817100,2.915365799) {$ M $};
\node[vertex,typeL,minimum size=4.4mm,font=\footnotesize] at (0.199817100,3.785365799) {$ L $};
\node[vertex,typeM,minimum size=4.4mm,font=\footnotesize] at (0.815000000,0.199817100) {$ M $};
\node[vertex,typeL,minimum size=4.4mm,font=\footnotesize] at (0.815000000,2.300182900) {$ L $};
\node[vertex,typeM,minimum size=4.4mm,font=\footnotesize] at (0.815000000,4.400548699) {$ M $};
\node[vertex,typeL,minimum size=4.4mm,font=\footnotesize] at (1.685000000,0.199817100) {$ L $};
\node[vertex,typeM,minimum size=4.4mm,font=\footnotesize] at (1.685000000,2.300182900) {$ M $};
\node[vertex,typeL,minimum size=4.4mm,font=\footnotesize] at (1.685000000,4.400548699) {$ L $};
\node[vertex,typeM,minimum size=4.4mm,font=\footnotesize] at (2.300182900,0.815000000) {$ M $};
\node[vertex,typeL,minimum size=4.4mm,font=\footnotesize] at (2.300182900,1.685000000) {$ L $};
\node[vertex,typeL,minimum size=4.4mm,font=\footnotesize] at (2.300182900,2.915365799) {$ L $};
\node[vertex,typeM,minimum size=4.4mm,font=\footnotesize] at (2.300182900,3.785365799) {$ M $};
\node[vertex,typeL,minimum size=4.4mm,font=\footnotesize] at (2.915365799,0.199817100) {$ L $};
\node[vertex,typeM,minimum size=4.4mm,font=\footnotesize] at (2.915365799,2.300182900) {$ M $};
\node[vertex,typeL,minimum size=4.4mm,font=\footnotesize] at (2.915365799,4.400548699) {$ L $};
\node[vertex,typeM,minimum size=4.4mm,font=\footnotesize] at (3.785365799,0.199817100) {$ M $};
\node[vertex,typeL,minimum size=4.4mm,font=\footnotesize] at (3.785365799,2.300182900) {$ L $};
\node[vertex,typeM,minimum size=4.4mm,font=\footnotesize] at (3.785365799,4.400548699) {$ M $};
\node[vertex,typeL,minimum size=4.4mm,font=\footnotesize] at (4.400548699,0.815000000) {$ L $};
\node[vertex,typeM,minimum size=4.4mm,font=\footnotesize] at (4.400548699,1.685000000) {$ M $};
\node[vertex,typeM,minimum size=4.4mm,font=\footnotesize] at (4.400548699,2.915365799) {$ M $};
\node[vertex,typeL,minimum size=4.4mm,font=\footnotesize] at (4.400548699,3.785365799) {$ L $};
\node[font=\footnotesize] at (2.300182900,2.300182900) {$4$};
\node[font=\footnotesize,text=black!60] at (1.250000000,1.250000000) {$8$};
\node[font=\footnotesize,text=black!60] at (3.350365799,1.250000000) {$8$};
\node[font=\footnotesize,text=black!60] at (1.250000000,3.350365799) {$8$};
\node[font=\footnotesize,text=black!60] at (3.350365799,3.350365799) {$8$};
\node[panel] at (2.300,0-.55) {(a) The coarsened $4.8.8$ map};
\draw[edge] (9.146751442,2.775000000) -- (8.475000000,3.446751442);
\draw[edge] (8.475000000,3.446751442) -- (7.525000000,3.446751442);
\draw[edge] (7.525000000,3.446751442) -- (6.853248558,2.775000000);
\draw[edge] (6.853248558,2.775000000) -- (6.853248558,1.825000000);
\draw[edge] (6.853248558,1.825000000) -- (7.525000000,1.153248558);
\draw[edge] (7.525000000,1.153248558) -- (8.475000000,1.153248558);
\draw[edge] (8.475000000,1.153248558) -- (9.146751442,1.825000000);
\draw[edge] (9.146751442,2.775000000) -- (9.146751442,1.825000000);
\draw[edge] (11.440254326,2.775000000) -- (10.768502884,3.446751442);
\draw[edge] (10.768502884,3.446751442) -- (9.818502884,3.446751442);
\draw[edge] (9.818502884,3.446751442) -- (9.146751442,2.775000000);
\draw[edge] (9.146751442,1.825000000) -- (9.818502884,1.153248558);
\draw[edge] (9.818502884,1.153248558) -- (10.768502884,1.153248558);
\draw[edge] (10.768502884,1.153248558) -- (11.440254326,1.825000000);
\draw[edge] (11.440254326,1.825000000) -- (11.440254326,2.775000000);
\draw[draw=black!35,dashed,line width=.7pt] (7.525000000,3.446751442) -- (8.475000000,1.153248558);
\draw[strongedge] (9.146751442,2.775000000) -- (6.853248558,1.825000000);
\draw[draw=black!35,dashed,line width=.7pt] (9.146751442,2.775000000) -- (11.440254326,1.825000000);
\draw[strongedge] (10.768502884,3.446751442) -- (9.818502884,1.153248558);
\node[vertex,typeL,minimum size=4.4mm,font=\footnotesize] at (6.853248558,1.825000000) {$ L $};
\node[vertex,typeM,minimum size=4.4mm,font=\footnotesize] at (6.853248558,2.775000000) {$ M $};
\node[vertex,typeM,minimum size=4.4mm,font=\footnotesize] at (7.525000000,1.153248558) {$ M $};
\node[vertex,typeL,minimum size=4.4mm,font=\footnotesize] at (7.525000000,3.446751442) {$ L $};
\node[vertex,typeL,minimum size=4.4mm,font=\footnotesize] at (8.475000000,1.153248558) {$ L $};
\node[vertex,typeM,minimum size=4.4mm,font=\footnotesize] at (8.475000000,3.446751442) {$ M $};
\node[vertex,typeM,minimum size=4.4mm,font=\footnotesize] at (9.146751442,1.825000000) {$ M $};
\node[vertex,typeL,minimum size=4.4mm,font=\footnotesize,draw=figureaccent,line width=1.1pt] at (9.146751442,2.775000000) {$ L $};
\node[vertex,typeL,minimum size=4.4mm,font=\footnotesize] at (9.818502884,1.153248558) {$ L $};
\node[vertex,typeM,minimum size=4.4mm,font=\footnotesize] at (9.818502884,3.446751442) {$ M $};
\node[vertex,typeM,minimum size=4.4mm,font=\footnotesize] at (10.768502884,1.153248558) {$ M $};
\node[vertex,typeL,minimum size=4.4mm,font=\footnotesize] at (10.768502884,3.446751442) {$ L $};
\node[vertex,typeL,minimum size=4.4mm,font=\footnotesize] at (11.440254326,1.825000000) {$ L $};
\node[vertex,typeM,minimum size=4.4mm,font=\footnotesize] at (11.440254326,2.775000000) {$ M $};
\draw[draw=figureaccent,line width=1pt] (9.146751442,2.775000000) -- (9.146751442,1.825000000);
\node[vertex,typeL,minimum size=4.4mm,font=\footnotesize,draw=figureaccent,line width=1.1pt] at (9.146751442,2.775000000) {$ L $};
\node[vertex,typeM,minimum size=4.4mm,font=\footnotesize] at (9.146751442,1.825000000) {$ M $};
\node[annotation] at (9.146751442,4.250000000) {exactly one octagon\\chooses this $L$};
\draw[arrow] (9.146751442,3.820000000) -- (9.146751442,3.075000000);
\draw[strongedge] (7.00,.70) -- (7.65,.70);
\node[font=\footnotesize,anchor=west] at (7.80,.70) {selected $LL$ diagonal};
\draw[draw=black!35,dashed,line width=.7pt] (7.00,.22) -- (7.65,.22);
\node[font=\footnotesize,anchor=west] at (7.80,.22) {unused alternative};
\node[panel] at (9.150,-.55) {(b) A diagonal choice forces its neighbor};
\end{tikzpicture}
\caption{Global propagation in the coarsened map. (a) Each vertex belongs to one quadrilateral and two octagons; regular polygons display the resulting combinatorics. (b) An alternating octagon has two pairs of opposite $L$ corners. Exactly one of the two octagons at a shared $L$ corner may choose that corner as an endpoint of its restored $LL$ diagonal. A choice in one octagon therefore forces the neighboring choice, and hence all choices in the connected octagon adjacency graph. Dashed diagonals show unused alternatives; they are not edges of the restored triangulation.}\label{fig:global-propagation}
\end{figure}

\begin{corollary}[Global identification]\label{cor:unique-packing-53}
Let $\PP$ be a triangulated disk packing of the whole plane with the words in \eqref{eq:ideal-words}. Suppose that its radii at vertices of types $S,M,L$ are $\lambda s,\lambda r,\lambda$, respectively, for one $\lambda>0$. Then $\PP$ is packing~53 up to a Euclidean similarity, possibly
reversing orientation.  In particular, it has two independent
translation periods.
\end{corollary}
\begin{proof}
Proposition~\ref{prop:unique-ideal-triangulation} identifies the
triangulation, including its vertex types, with that of packing~53.
Scale packing~53 by $\lambda$. Corresponding tangency edges then have equal
lengths, since an edge length is the sum of its endpoint radii.
Place one corresponding triangular face by a Euclidean isometry.
Across any edge, the next triangle is uniquely determined by its
three side lengths and its position on the opposite side of that edge.
Connectedness of the face adjacency graph propagates this placement
to every face.  Thus the disk centers and radii agree throughout
the plane after one similarity.  Packing~53 has two independent
translation periods, and this property is preserved by similarities.
\end{proof}

\subsection{Proofs of the main theorems}
\begin{proof}[Proof of Theorem~\ref{thm:structure}]
Normalize $\sup_vr_v=1$, so every actual radius lies in $[q_0,1]$. Proposition~\ref{prop:CZ} gives degrees $5,6,7$, and Theorem~\ref{thm:no-exceptional} excludes every exceptional corona. Set
\[
h_v=\log(r_v/\rho_v),
\]
where $\boldsymbol\rho$ is \eqref{eq:reference-assignment}. The actual and reference radii belong to a common compact positive interval, so $h$ is bounded. Lemma~\ref{lem:jacobian} and Proposition~\ref{prop:reference-excess} give
\begin{equation}\label{eq:final-laplacian}
\Delta_w h(v)=2\pi-A_v(\boldsymbol\rho)\le0.
\end{equation}
The weights are positive, symmetric, and uniformly bounded above. Lemma~\ref{lem:recurrence} therefore makes $h$ constant. There is a single $\lambda>0$ with $r_v=\lambda\rho_v$ for every $v$.

Lemma~\ref{lem:degree-five} supplies an $S$ vertex. Its word contains both $M$ and $L$, so all three reference radii occur. This proves \eqref{eq:forced-radii} and $q(\PP)=s$. In the normalization used here, in fact $\lambda=1$.

Equation~\eqref{eq:final-laplacian} now has zero left side, so every word is an equality word from Proposition~\ref{prop:reference-excess}. We only have to exclude passive $M$ vertices. A passive vertex has only $M$ neighbors. An active $M$ has a unique $M$ neighbor, flanked by $L$; the two incident triangles imply that this neighbor has $L$ neighbors and is not passive. Thus passive vertices have no edge to an active $M$, an $S$, or an $L$. Since $G$ is connected and contains an $S$, no passive vertex occurs. The words are exactly those in \eqref{eq:ideal-words}.

Corollary~\ref{cor:unique-packing-53} now identifies $\PP$ with packing~53 up to similarity. The latter has two independent translation periods, so the same is true of $\PP$. This completes the proof.
\end{proof}

\begin{proof}[Proof of Theorem~\ref{thm:optimal}]
If $q(\PP)<q_0$, then $q(\PP)<s$ by \eqref{eq:threshold-order}. Otherwise Theorem~\ref{thm:structure} identifies $\PP$ with packing~53 up to similarity and gives $q(\PP)=s$. Packing~53 supplies an example attaining $s$, and its similarities are the only equality cases.
\end{proof}

\section{Further questions}\label{sec:questions}
The preceding argument identifies every nonuniform triangulated packing with ratio at least $q_0$, including all equality packings, with packing~53 up to similarity. Thus their global triangulation and double periodicity are determined.

A remaining question concerns the gap below $s$. The structure theorem holds at the certified cutoff $q_0$. To lower this cutoff within the present proof, all local exclusions and comparison hypotheses must remain valid at the new value.

\begin{problem}
Determine the supremum of the radius ratios of nonuniform triangulated packings whose ratio is different from $s$.
\end{problem}

There is also a broader question about density when the triangulation
assumption is dropped. In the following discussion, $\PP$ may be any disk
packing satisfying \eqref{eq:bounded-radii}. Write
\begin{equation}\label{eq:general-density}
\underline\delta(\PP):=
\liminf_{R\to\infty}
\frac{\operatorname{area}\bigl(\overline B(0,R)\cap
\bigcup_{D\in\PP}D\bigr)}{\pi R^2},
\qquad
\delta_{\mathrm{hex}}:=\frac{\pi}{\sqrt{12}}.
\end{equation}
If one asks only for $\underline\delta(\PP)\ge\delta_{\mathrm{hex}}$,
the maximum ratio is $1$, attained by the hexagonal packing. Even among
nonuniform packings the supremum is $1$: shrinking one disk of the
unit-radius hexagonal packing to any radius $t\in(0,1)$ gives ratio $t$
and changes the covered area by only the fixed amount $\pi(1-t^2)$.
Its limiting density is therefore unchanged. The substantive question
requires a strict improvement in density; see
Fernique~\cite[Section~4, Open Question~4]{fernique-survey}.

\begin{problem}[Density without triangulation]\label{prob:density-threshold}
Determine
\begin{equation}\label{eq:density-threshold}
q_{\mathrm{dens}}:=
\sup\bigl\{q(\PP):\PP\text{ satisfies \eqref{eq:bounded-radii} and }
\underline\delta(\PP)>\delta_{\mathrm{hex}}\bigr\},
\end{equation}
where no triangulation, periodicity, or restriction on the number of
radii is imposed.
\end{problem}

The use of a supremum is necessary. Suppose a packing has lower density
$d>\delta_{\mathrm{hex}}$, and put $a=\inf_v r_v$ and $b=\sup_v r_v$.
Here $a<b$, since congruent disks have density at most
$\delta_{\mathrm{hex}}$.
Choose $c\in(a,b)$ sufficiently close to $b$ that
$(c/b)^2d>\delta_{\mathrm{hex}}$, and replace each radius $r_v$ by
$\min\{r_v,c\}$, keeping all centers fixed. The resulting packing
$\PP'$ has ratio $a/c>a/b$ and
\[
\underline\delta(\PP')\ge(c/b)^2\underline\delta(\PP)
>\delta_{\mathrm{hex}}.
\]
Indeed, every disk retains at least the fraction $(c/b)^2$ of its area,
and disks meeting the boundary of $\overline B(0,R)$ contribute only
$O(R)$ to the covered-area comparison because the radii are bounded.
Thus no packing with strictly improved density attains $q_{\mathrm{dens}}$.

Connelly and Pierre~\cite[Section~4, Theorem~4.1]{connelly-pierre}
give a doubly periodic deformation of packing~53. Let
$\delta_{\mathrm{CP}}(t)$ denote the density of this family when its
radius ratio is $t$, as defined by equations~(4.1)--(4.2) of that paper.
Let $q_{\mathrm{CP}}$ be the first parameter $t>s$ on this branch
for which $\delta_{\mathrm{CP}}(t)=\delta_{\mathrm{hex}}$.
Numerically,
\[
q_{\mathrm{CP}}=0.6585340820\ldots.
\]
The density is strictly greater than $\delta_{\mathrm{hex}}$ for
$s\le t<q_{\mathrm{CP}}$. Hence this construction gives the lower
bound $q_{\mathrm{dens}}\ge q_{\mathrm{CP}}$, even though the
packing at $q_{\mathrm{CP}}$ has density exactly $\delta_{\mathrm{hex}}$.

For an upper bound, the finite-packing theorem of Blind and
Blind~\cite[Theorem~1]{blind-blind} implies that a packing with
ratio at least
\begin{equation}\label{eq:blind-threshold}
q_B:=\sqrt{\frac{\sqrt{12}-7\tan(\pi/7)}
{5\tan(\pi/5)-\sqrt{12}}}
=0.7429909632\ldots
\end{equation}
has lower density at most $\delta_{\mathrm{hex}}$.
To pass from their finite statement to the whole plane, apply it to
the disks contained in $\overline B(0,R)$ and let $R\to\infty$;
the omitted boundary disks contribute an area of order $R$.
The bound applies to arbitrarily many radii. Consequently,
\begin{equation}\label{eq:density-threshold-bounds}
0.6585340820\ldots=q_{\mathrm{CP}}
\le q_{\mathrm{dens}}\le q_B=0.7429909632\ldots.
\end{equation}
Connelly and Pierre~\cite[Conjecture~4.2]{connelly-pierre}
conjecture that the lower bound is sharp. Theorem~\ref{thm:optimal}
settles the triangulated assertion in that conjecture; determining
$q_{\mathrm{dens}}$ additionally requires control of packings whose
tangency graphs are not triangulations.

A related question fixes the permitted interval of radii and asks for
the greatest possible density. Let $\delta_{53}$ denote the density of
packing~53 with the exact reference radii $s,r,1$ from
\eqref{eq:reference-numbers}. Connelly and
Pierre~\cite[Section~4]{connelly-pierre} give
\[
\delta_{53}=\delta_{\mathrm{CP}}(s)=0.909306501624174\ldots.
\]

\begin{problem}[Density on the optimal radius interval]
\label{prob:interval-density}
Does every disk packing $\PP$ whose radii satisfy $s\le r_v\le1$ satisfy
\[
\underline\delta(\PP)\le\delta_{53}?
\]
No triangulation, periodicity, or restriction on the number of distinct
radii is assumed.
\end{problem}

Fernique~\cite[Theorem~1]{fernique2019} proves this inequality when
every radius belongs to the fixed set $\{s,r,1\}$.
Problem~\ref{prob:interval-density} allows every radius in $[s,1]$.
Together with the uniform case, Theorem~\ref{thm:structure} settles the
question within the triangulated class. Thus a packing answering Problem~\ref{prob:interval-density}
negatively would have to be nontriangulated and use a radius outside
$\{s,r,1\}$.

Florian's bound, stated in Connelly and
Pierre~\cite[Theorem~1.1]{connelly-pierre}, gives, for every such packing,
\begin{equation}\label{eq:interval-density-bound}
\underline\delta(\PP)\le
\frac{\pi s^2+2(1-s^2)\arcsin\!\left(\dfrac{s}{1+s}\right)}
{2s\sqrt{1+2s}}
=0.910425244015304\ldots.
\end{equation}
This upper bound exceeds $\delta_{53}$ by approximately $0.0011187424$.
The Connelly--Pierre deformation does not improve on $\delta_{53}$:
its density $\delta_{\mathrm{CP}}(t)$ strictly decreases along the branch
$s\le t\le q_{\mathrm{CP}}$.

Finally, the exceptional-corona certificate uses only completed-flower angle equations on an unfolded disk. It would be useful to replace the large interval trace by a shorter analytic inequality or a substantially smaller certificate. Such a reduction would simplify the verification without changing the geometric obstruction.

\section*{Disclosure of AI assistance}
The author developed the proof strategy in 2019; its finite-subgraph exclusion approach was subsequently outlined in joint work with Robert Connelly~\cite{cz}. At that time, the author was unable to implement the necessary interval-arithmetic computations efficiently enough to complete the proof. The present paper is the result of substantial human--AI collaboration. The author supplied the underlying mathematical strategy and geometric arguments, while ChatGPT (OpenAI) produced efficient implementations of the interval-arithmetic searches and exact verification programs that made the finite exclusions computationally feasible. Both contributions were essential to completing the proof in its present form. After the proof outline was established, ChatGPT was also used to draft and revise parts of the text and to generate figures. The author takes full responsibility for the mathematical claims, computational results, and final presentation.

\section*{Data availability}
To the best of the author's knowledge, all data presented in this paper are included
in the paper itself or in its accompanying supporting file package. Additional data
can be made available by the author upon reasonable request.

\appendix
\section{The degree-five and degree-six exclusions}\label{app:local}
This appendix gives the finite configurations for Proposition~\ref{prop:local}. Every search here uses the exact lower radius bound $q_0=65105/100000$. The interval rules and their soundness are established in Appendix~\ref{app:arithmetic}.

\subsection{Geometric coverage of the configurations}\label{app:patches}
\begin{figure}[htbp]
\centering
\begin{tikzpicture}[font=\small]
\draw[edge] (0.00000000,0.00000000) -- (0.00000000,1.25000000);
\draw[edge] (0.00000000,1.25000000) -- (-1.18882065,0.38627124);
\draw[edge] (0.00000000,1.25000000) -- (-1.31663897,1.81219807);
\draw[edge] (0.00000000,1.25000000) -- (1.31663897,1.81219807);
\draw[edge] (-1.31663897,1.81219807) -- (-0.53433305,2.51383933);
\draw[edge] (-0.53433305,2.51383933) -- (0.53433305,2.51383933);
\draw[edge] (0.53433305,2.51383933) -- (1.31663897,1.81219807);
\draw[edge] (0.00000000,1.25000000) -- (-0.53433305,2.51383933);
\draw[edge] (0.00000000,1.25000000) -- (0.53433305,2.51383933);
\draw[edge] (0.00000000,0.00000000) -- (-1.18882065,0.38627124);
\draw[edge] (-1.18882065,0.38627124) -- (-0.73473157,-1.01127124);
\draw[edge] (-1.18882065,0.38627124) -- (-2.13036660,-0.69219807);
\draw[edge] (-1.18882065,0.38627124) -- (-1.31663897,1.81219807);
\draw[edge] (-2.13036660,-0.69219807) -- (-2.44421525,0.79417368);
\draw[edge] (-2.44421525,0.79417368) -- (-1.31663897,1.81219807);
\draw[edge] (-1.18882065,0.38627124) -- (-2.44421525,0.79417368);
\draw[edge] (0.00000000,0.00000000) -- (-0.73473157,-1.01127124);
\draw[edge] (-0.73473157,-1.01127124) -- (0.73473157,-1.01127124);
\draw[edge] (-0.73473157,-1.01127124) -- (-0.00000000,-2.24000000);
\draw[edge] (-0.73473157,-1.01127124) -- (-2.13036660,-0.69219807);
\draw[edge] (-0.00000000,-2.24000000) -- (-1.04531317,-2.34781183);
\draw[edge] (-1.04531317,-2.34781183) -- (-1.90988220,-1.71966566);
\draw[edge] (-1.90988220,-1.71966566) -- (-2.13036660,-0.69219807);
\draw[edge] (-0.73473157,-1.01127124) -- (-1.04531317,-2.34781183);
\draw[edge] (-0.73473157,-1.01127124) -- (-1.90988220,-1.71966566);
\draw[edge] (0.00000000,0.00000000) -- (0.73473157,-1.01127124);
\draw[edge] (0.73473157,-1.01127124) -- (1.18882065,0.38627124);
\draw[edge] (0.73473157,-1.01127124) -- (2.13036660,-0.69219807);
\draw[edge] (0.73473157,-1.01127124) -- (-0.00000000,-2.24000000);
\draw[edge] (2.13036660,-0.69219807) -- (1.90988220,-1.71966566);
\draw[edge] (1.90988220,-1.71966566) -- (1.04531317,-2.34781183);
\draw[edge] (1.04531317,-2.34781183) -- (-0.00000000,-2.24000000);
\draw[edge] (0.73473157,-1.01127124) -- (1.90988220,-1.71966566);
\draw[edge] (0.73473157,-1.01127124) -- (1.04531317,-2.34781183);
\draw[edge] (0.00000000,0.00000000) -- (1.18882065,0.38627124);
\draw[edge] (1.18882065,0.38627124) -- (0.00000000,1.25000000);
\draw[edge] (1.18882065,0.38627124) -- (1.31663897,1.81219807);
\draw[edge] (1.18882065,0.38627124) -- (2.13036660,-0.69219807);
\draw[edge] (1.31663897,1.81219807) -- (2.44421525,0.79417368);
\draw[edge] (2.44421525,0.79417368) -- (2.13036660,-0.69219807);
\draw[edge] (1.18882065,0.38627124) -- (2.44421525,0.79417368);
\node[vertex,minimum size=6.8mm,typeL] at (0.00000000,1.25000000) {$v_0$};
\node[vertex,minimum size=6.8mm,typeM] at (-1.18882065,0.38627124) {$v_1$};
\node[vertex,minimum size=6.8mm,typeL] at (-0.73473157,-1.01127124) {$v_2$};
\node[vertex,minimum size=6.8mm,typeL] at (0.73473157,-1.01127124) {$v_3$};
\node[vertex,minimum size=6.8mm,typeM] at (1.18882065,0.38627124) {$v_4$};
\node[vertex,minimum size=6.8mm] at (-1.31663897,1.81219807) {$w_0$};
\node[vertex,minimum size=6.8mm] at (-2.13036660,-0.69219807) {$w_1$};
\node[vertex,minimum size=6.8mm] at (-0.00000000,-2.24000000) {$w_2$};
\node[vertex,minimum size=6.8mm] at (2.13036660,-0.69219807) {$w_3$};
\node[vertex,minimum size=6.8mm] at (1.31663897,1.81219807) {$w_4$};
\node[circle,fill=black,inner sep=0pt,minimum size=2.1mm] at (-0.53433305,2.51383933) {};
\node[font=\footnotesize] at (-0.62165596,2.92466133) {$e_{0,1}$};
\node[circle,fill=black,inner sep=0pt,minimum size=2.1mm] at (0.53433305,2.51383933) {};
\node[font=\footnotesize] at (0.62165596,2.92466133) {$e_{0,2}$};
\node[circle,fill=black,inner sep=0pt,minimum size=2.1mm] at (-2.44421525,0.79417368) {};
\node[font=\footnotesize] at (-2.84365898,0.92396081) {$e_{1,1}$};
\node[circle,fill=black,inner sep=0pt,minimum size=2.1mm] at (-1.04531317,-2.34781183) {};
\node[font=\footnotesize] at (-1.21614256,-2.73150092) {$e_{2,1}$};
\node[circle,fill=black,inner sep=0pt,minimum size=2.1mm] at (-1.90988220,-1.71966566) {};
\node[font=\footnotesize] at (-2.22200303,-2.00070051) {$e_{2,2}$};
\node[circle,fill=black,inner sep=0pt,minimum size=2.1mm] at (1.90988220,-1.71966566) {};
\node[font=\footnotesize] at (2.22200303,-2.00070051) {$e_{3,1}$};
\node[circle,fill=black,inner sep=0pt,minimum size=2.1mm] at (1.04531317,-2.34781183) {};
\node[font=\footnotesize] at (1.21614256,-2.73150092) {$e_{3,2}$};
\node[circle,fill=black,inner sep=0pt,minimum size=2.1mm] at (2.44421525,0.79417368) {};
\node[font=\footnotesize] at (2.84365898,0.92396081) {$e_{4,1}$};
\node[vertex,typeS,minimum size=6.8mm] at (0.00000000,0.00000000) {$A$};
\node[panel] at (0.00000000,-3.52000000) {(a) A labeled wheel and its first extensions};
\draw[edge] (7.25000000,0.35000000) -- (5.15000000,1.85000000);
\draw[edge] (7.25000000,0.35000000) -- (6.43000000,2.55000000);
\draw[edge] (7.25000000,0.35000000) -- (8.07000000,2.55000000);
\draw[edge] (7.25000000,0.35000000) -- (9.35000000,1.85000000);
\draw[edge] (5.15000000,1.85000000) -- (6.43000000,2.55000000);
\draw[edge] (6.43000000,2.55000000) -- (8.07000000,2.55000000);
\draw[edge] (8.07000000,2.55000000) -- (9.35000000,1.85000000);
\node[font=\footnotesize,fill=white,draw=black!70,rounded corners=1pt,inner sep=3pt] at (5.15000000,1.85000000) {$e_{i,1}$};
\node[font=\footnotesize,fill=white,draw=black!70,rounded corners=1pt,inner sep=3pt] at (6.43000000,2.55000000) {$v_i$};
\node[font=\footnotesize,fill=white,draw=black!70,rounded corners=1pt,inner sep=3pt] at (8.07000000,2.55000000) {$v_{i+1}$};
\node[font=\footnotesize,fill=white,draw=black!70,rounded corners=1pt,inner sep=3pt] at (9.35000000,1.85000000) {$e_{i+1,\mathrm{last}}$};
\node[vertex,minimum size=8mm] at (7.25000000,0.35000000) {$w_i$};
\draw[draw=black!45,densely dashed] (5.15000000,1.85000000) .. controls (4.25000000,-1.05000000) and (10.25000000,-1.05000000) .. (9.35000000,1.85000000);
\node[annotation] at (7.25000000,-0.87000000) {unknown remainder of the flower};
\node[panel] at (7.25000000,-3.52000000) {(b) A known consecutive path};
\end{tikzpicture}
\caption{The local incidence data used in Appendix~\ref{app:patches}. Panel (a) displays one of the sixteen degree assignments, with $(d(v_0),\ldots,d(v_4))=(7,6,7,7,6)$; the extra vertices complete each $v_i$-flower. Panel (b) isolates the known four-entry path at $w_i$. These are labeled combinatorial diagrams. Exterior labels are kept independent in the exclusion problem even when some of them could denote the same disk in a global packing; the drawing does not impose a simple outer boundary.}\label{fig:verification-wheel}
\end{figure}
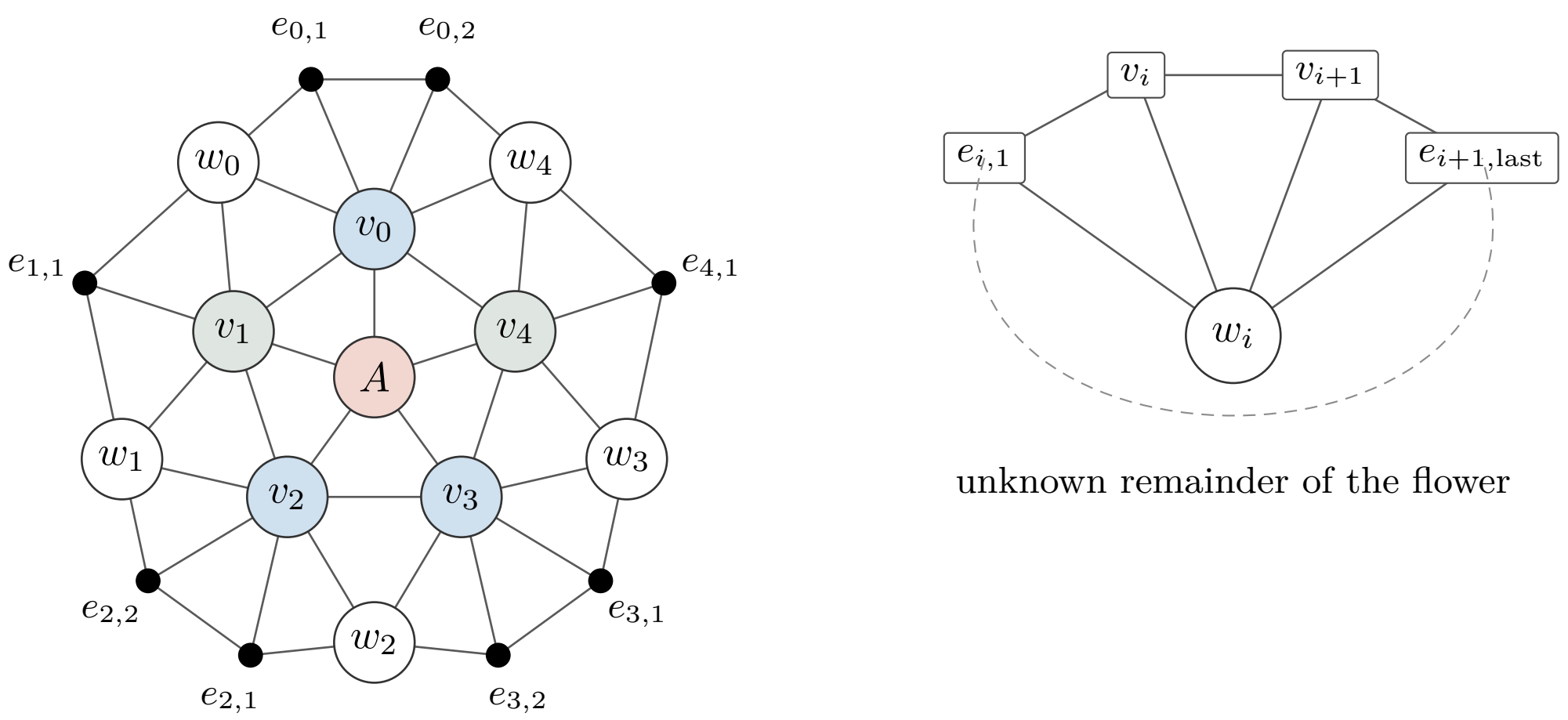

Choose a degree-five center $A$ with cyclic neighbors $v_0,\ldots,v_4$. Proposition~\ref{prop:CZ} allows us to choose $v_0$ of degree seven; the other neighbors have degree six or seven. Subscripts in this paragraph are taken modulo five. Let $w_i$ be the vertex opposite $A$ across $v_iv_{i+1}$. The complete list at $v_i$ is
\begin{equation}\label{eq:wheel-list}
[w_{i-1},v_{i-1},A,v_{i+1},w_i,e_{i,1},\ldots,e_{i,d(v_i)-5}].
\end{equation}
There are one or two extra labels according as $d(v_i)=6$ or $7$. The known consecutive path in the flower at $w_i$ is
\begin{equation}\label{eq:opposite-path}
[e_{i,1},v_i,v_{i+1},e_{i+1,d(v_{i+1})-5}].
\end{equation}
Each extra vertex has the three-entry path consisting of its predecessor in the outer chain, $v_i$, and its successor. In Figure~\ref{fig:verification-wheel}, $e_{j,\mathrm{last}}$ abbreviates $e_{j,d(v_j)-5}$.

The data are local face incidences, not a claim that the complete labeled outer boundary embeds as a simple polygon. Exterior labels may have the same image in a packing. Their radii and degrees are allowed to vary independently in the interval problem; any actual identification only adds equality constraints that we do not need.

Within a prescribed consecutive link path the actual neighbors are distinct. Consecutive faces move in one direction around the target link, as in Lemma~\ref{lem:link}; a return before the allowed degree would close the link too early. In the extended problems a partial path has at most five vertices. A path with five vertices may be closed when the tested target degree is five. At every two-face edge used for Proposition~\ref{prop:CZ}(iii), a completed endpoint flower identifies the two distinct opposite face vertices.

To complete a partial flower, append the missing neighbors, extend the two endpoint paths on the corresponding sides, and retain the predecessor--center--successor path at each fresh label. Every constraint used by the searches comes from these face incidences, from Proposition~\ref{prop:CZ}, or from a local statement already proved. No unspecified global non-overlap condition is imposed.

For the sixteen base cases, choose bits $b_0,b_1,b_2,b_3\in\{0,1\}$ and set
\begin{equation}\label{eq:mask}
d(v_0)=7,\qquad d(v_i)=6+b_{i-1}\ (1\le i\le4),
\qquad m=\sum_{j=0}^3 2^j b_j.
\end{equation}
When the orientation is otherwise free, reflection allows the radius ordering $r_{v_1}\le r_{v_4}$.

\subsection{Forcing the degree-five word and partner}\label{app:pairing-cases}
For each degree assignment in \eqref{eq:mask}, we test the angle equations at the six completed flowers $A,v_0,\ldots,v_4$ and necessary angle bounds at the partial exterior flowers. We also impose Proposition~\ref{prop:CZ}, $r_A\le0.702$, and the orientation convention $r_{v_1}\le r_{v_4}$. The exterior radius variables remain independent.

\Needspace{7\baselineskip}
A \emph{box} is a product $B=\prod_x[\ell_x,u_x]$ of radius intervals. Monotonicity of $T$ bounds each required angle sum on $B$, including possible completions of partial flowers. A box is rejected only by a certified contradiction with a radius, angle, or degree condition; otherwise its radius or degree bounds are tightened, or a radius interval is split into two closed children. The exact comparisons and covering argument are given in Appendix~\ref{app:arithmetic}.

The searches exclude the following twelve cases; the table gives the number of boxes processed:
\begin{center}
\begin{tabular}{@{}rrrrrrrr@{}}
\toprule
$m$ & Boxes & $m$ & Boxes & $m$ & Boxes & $m$ & Boxes\\
\midrule
0&1&3&1&7&1&13&515\\
1&1&4&107&8&75&14&603\\
2&139&5&1&11&511&15&13\\
\bottomrule
\end{tabular}
\end{center}
Every search ends without a surviving or pending box; a count of one requires no split. These are rooted, oriented cases subject to $r_{v_1}\le r_{v_4}$, so different cases can represent the same cyclic degree word.

Of the four remaining cases, $10$ can be re-rooted at its isolated degree-seven neighbor to give case $6$, and $12$ at the middle of its three consecutive degree-seven neighbors to give case $9$. Both target words are invariant under reversal, so the radius-order convention remains available.

Completing the flower at $w_0$ in case $9$ excludes its possible degrees $5,6,7$ in $623,195,1$ boxes, respectively. Completing the flower at $w_2$ in case $6$ excludes degrees $6,7$ in $1067,3$ boxes; hence $d(w_2)=5$. No box remains unresolved.

Thus $A$ has cyclic degree word $7,6,7,7,6$, and the vertex opposite $A$ across its consecutive $7$--$7$ edge also has degree five. This proves Proposition~\ref{prop:local}(i), and the involution in Corollary~\ref{cor:pairing} follows.

\subsection{Excluding three degree-five neighbors}\label{app:three-small}
Suppose a degree-six vertex has three degree-five neighbors. No two $S$ are adjacent, and each $S$ is flanked by $L$, so its word is $\word{LSLSLS}$. Root at one of the three $S$ vertices, choose its isolated $L$ as $v_0$, and reflect so that the proposed $M$ is $v_1$. This fixes the orientation. In this test the radius-order inequality $r_{v_1}\le r_{v_4}$ is \emph{not} imposed.

To specify the labels used in the certificate, set $A=0$, $v_i=1+i$, and $w_i=6+i$. In case $6$, assign the extra labels successively from $11$ onward. The paired $S$ is vertex $8$, and completing it introduces vertex $19$, with list
\[
[14,3,4,17,19].
\]
Vertices $3,4$ form its unique consecutive $LL$ pair. The remaining degrees at $14,17,19$ are therefore $6,6,7$. The proposed $M$ has list
\[
[6,1,0,3,7,13].
\]
Alternation forces $6,7$ to be degree five and $13$ to be degree seven. Complete the flowers at $8,14,17,19,6,7$ with degrees $5,6,6,7,5,5$, respectively. Propagate Proposition~\ref{prop:local}(i) where its required faces are known.

The resulting interval problem is excluded in five boxes. The computation verifies $467$ strict angle-sum comparisons and $1174$ distinct angular enclosures. Its minimum strict margin beyond the enclosing bound for $2\pi$ is $13\cdot10^{-6}$ radians. Together with the covering principle of Appendix~\ref{app:arithmetic}, this proves Proposition~\ref{prop:local}(ii).

\section{Exact interval arithmetic}\label{app:arithmetic}
We give the proof rules used by the local certificates and by the new exceptional-corona certificate. All interval endpoints represent bounds on real variables. They do not restrict the radii to a discrete grid.

\subsection{Monotone bounds and complete coverage}
For a finite set of $N$ radius variables, let
\[
B=\prod_{j=0}^{N-1}[\ell_j,u_j]
\]
be a closed interval box. For a completed center $i$ with neighbor list $J_i=(j_1,\ldots,j_{d_i})$, set
\[
F_i(a;b_1,\ldots,b_{d_i})=\sum_{k=1}^{d_i}T(a,b_k,b_{k+1}),
\qquad b_{d_i+1}=b_1.
\]
Write $\underline{\mathbf b}_i=(\ell_{j_1},\ldots,\ell_{j_{d_i}})$ and $\overline{\mathbf b}_i=(u_{j_1},\ldots,u_{j_{d_i}})$. Monotonicity of $T$ gives
\begin{equation}\label{eq:box-bounds}
F_i(u_i;\underline{\mathbf b}_i)
\le\calA_i(\mathbf a)
\le F_i(\ell_i;\overline{\mathbf b}_i)
\qquad(\mathbf a\in B).
\end{equation}
Here $\calA_i$ denotes the sum for that completed flower, as in \eqref{eq:finite-angle}. If the lower bound exceeds $2\pi$ or the upper bound is less than $2\pi$, the box contains no solution.

The same monotonicity allows contractions. For a neighbor label $j\in J_i$, the notation $\mathbf b[j\leftarrow t]$ means that its entry in the neighbor list is replaced by $t$. The four rules are
\begin{center}
\small
\begin{tabular}{@{}p{.67\linewidth}p{.25\linewidth}@{}}
\toprule
Verified strict inequality & Necessary consequence\\
\midrule
$F_i(t;\underline{\mathbf b}_i)>2\pi$ & $a_i>t$\\[3pt]
$F_i(t;\overline{\mathbf b}_i)<2\pi$ & $a_i<t$\\[3pt]
$F_i(\ell_i;\overline{\mathbf b}_i[j\leftarrow t])<2\pi$ & $a_j>t$\\[3pt]
$F_i(u_i;\underline{\mathbf b}_i[j\leftarrow t])>2\pi$ & $a_j<t$\\
\bottomrule
\end{tabular}
\end{center}
Keeping $t$ as a closed endpoint after a strict contraction is harmless: it only enlarges the retained set.

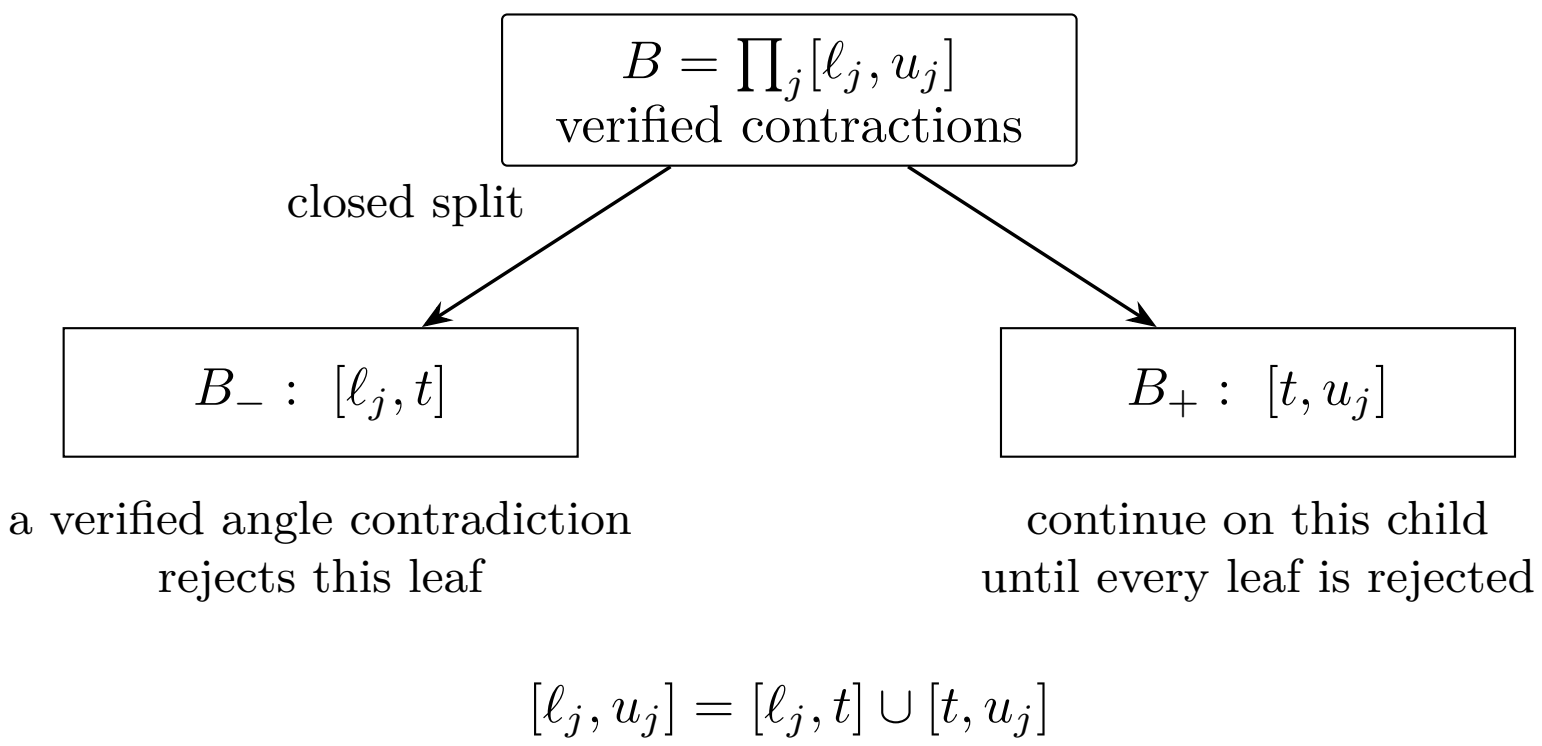
\begin{figure}[htbp]
\centering
\begin{tikzpicture}[font=\small]
\node[draw,rounded corners=1pt,minimum width=3.8cm,minimum height=1.0cm,align=center] (a) at (0,0) {$B=\prod_j[\ell_j,u_j]$\\verified contractions};
\node[draw,minimum width=3.4cm,minimum height=.85cm] (b) at (-3.1,-2) {$B_-:\ [\ell_j,t]$};
\node[draw,minimum width=3.4cm,minimum height=.85cm] (c) at (3.1,-2) {$B_+:\ [t,u_j]$};
\draw[arrow] (a)--node[above left,font=\footnotesize] {closed split} (b);
\draw[arrow] (a)--(c);
\node[align=center,font=\footnotesize] at (-3.1,-3.05) {a verified angle contradiction\\rejects this leaf};
\node[align=center,font=\footnotesize] at (3.1,-3.05) {continue on this child\\until every leaf is rejected};
\node[align=center] at (0,-4.1) {$[\ell_j,u_j]=[\ell_j,t]\cup[t,u_j]$};
\end{tikzpicture}
\caption{The covering principle for the finite interval proof. A solution survives every valid contraction and belongs to at least one child of each closed split. Consequently, it cannot be lost at a shared endpoint. The drawing is schematic; the recorded exceptional-corona certificate has $3894$ splits and $3895$ rejected leaves.}\label{fig:interval-cover}
\end{figure}

A binary split replaces one interval $[\ell_j,u_j]$ by both $[\ell_j,t]$ and $[t,u_j]$, with $\ell_j<t<u_j$. The children cover the parent, including their shared endpoint; see Figure~\ref{fig:interval-cover}.

\begin{lemma}[Certificate covering principle]\label{lem:certificate-soundness}
Suppose a finite tree starts from the full prescribed box and uses only valid contractions, closed splits, and rejections of the form \eqref{eq:box-bounds}. If every leaf is rejected, the system of completed-flower equations has no solution in the initial box.
\end{lemma}
\begin{proof}
A solution in a current box survives every valid contraction and belongs to at least one child of a split. It cannot belong to a rejected box, where one angle sum is proved different from $2\pi$. Following a putative solution down the finite tree would therefore give an unrejected leaf, a contradiction.
\end{proof}

\subsection{Partial flowers and degree choices}
For the older local exclusions, some vertices have only a consecutive neighbor path and an allowed degree set contained in $\{5,6,7\}$. Given a tested degree, complete its list by unknown placeholders. For a lower angle bound give each placeholder the minimum allowed radius; for an upper bound give it the maximum. Known radii and the center are replaced by the appropriate interval endpoints. Every actual completion is included in these bounds, whether or not some placeholders coincide with labels elsewhere.

The programs also compare the smallest and largest permitted degrees. Inserting an additional placeholder increases either extremal sum. If a placeholder is already present, the insertion adds a positive angle between two identical endpoint radii. If none is present, one angle is replaced by two; the exact inequality
\begin{equation}\label{eq:insertion}
2T(1,q_0,q_0)>1.621330>1.301230>T(q_0,1,1)
\end{equation}
ensures a positive increase in this case as well. Thus taking the smallest degree for a lower bound and the largest for an upper bound is valid. An impossible degree is removed only after a necessary angle comparison or an already proved degree restriction excludes it.

Every degree-five center has radius less than $0.702$, since
\[
5T(0.702,1,1)<2\pi.
\]
This comparison is checked exactly before the local searches. The prohibitions on adjacent degree-five vertices and on four degree-seven vertices across a two-face edge are propagated using Proposition~\ref{prop:CZ}. The degree-five rule is used in the three-small-neighbor exclusion only after the pairing certificate has established it.

Radius endpoints in the supplied traces are integer multiples of $10^{-5}$. If a box reaches the minimum implemented width without an actual contradiction, it must be reported as surviving, not rejected. The recorded local searches have no such surviving boxes. This observation is important: finite precision in the endpoint representation is not an inference of infeasibility.

\subsection{Rational angle enclosures}
The local generating programs certify angle bounds by sine polynomials. The separate exceptional-corona replay uses different cosine polynomials. We record both elementary tests.

\begin{lemma}[Sine enclosures]\label{lem:sine-enclosure}
For $0<x<1$, put
\[
S_9(x)=\sum_{k=0}^4\frac{(-1)^kx^{2k+1}}{(2k+1)!},
\qquad
S_{11}(x)=\sum_{k=0}^5\frac{(-1)^kx^{2k+1}}{(2k+1)!}.
\]
Then $S_{11}(x)<\sin x<S_9(x)$. For positive $a,b,c$,
\begin{enumerate}
\item $S_9(x)>0$ and $S_9(x)^2(a+b)(a+c)<bc$ imply $2x<T(a,b,c)$;
\item $S_{11}(x)>0$ and $S_{11}(x)^2(a+b)(a+c)>bc$ imply $T(a,b,c)<2x$.
\end{enumerate}
\end{lemma}
\begin{proof}
The alternating sine series gives the strict polynomial bounds. Apply the half-angle formula \eqref{eq:angle} and strict monotonicity of sine on $(0,\pi/2)$.
\end{proof}

To certify a full-angle endpoint $m/10^6$, the sine-based programs take $x=m/(2\cdot10^6)$ and check $0<x<1$. All denominators are cleared, leaving exact integer comparisons. Floating-point trigonometry only proposes endpoints and comparisons to try; every effective contraction or rejection must pass the integer test.

\begin{lemma}[Cosine enclosures]\label{lem:cosine-enclosure}
For $n\ge0$, define
\[
C_{2n}(x)=\sum_{k=0}^n\frac{(-1)^kx^{2k}}{(2k)!}.
\]
For $0<x<2$, one has $C_{14}(x)<\cos x<C_{12}(x)$. Put
\[
\Gamma(a,b,c)=1-\frac{2bc}{(a+b)(a+c)}=\cos T(a,b,c).
\]
Then
\begin{align*}
C_{14}(x)>\Gamma(a,b,c)&\ \Longrightarrow\ x<T(a,b,c),\\
C_{12}(x)<\Gamma(a,b,c)&\ \Longrightarrow\ T(a,b,c)<x.
\end{align*}
For rational arguments, these tests can be verified by integer arithmetic.
\end{lemma}
\begin{proof}
In each cosine tail, the absolute values of the terms decrease and their signs alternate; their ratio is at most $4/((2k+1)(2k+2))$ beyond the first omitted term. This proves the strict polynomial enclosures. Both $x$ and $T(a,b,c)$ lie in $(0,\pi)$, where cosine is strictly decreasing. Comparing their cosines gives the two implications.
\end{proof}

For the exact implementation, let $A,B,C$ be integer radius endpoints in units of $10^{-5}$ and let $x=m/10^6$, with $0<m<2\cdot10^6$. Define
\[
H=(A+B)(A+C),\qquad N=H-2BC,\qquad D=14!\,(10^6)^{14}.
\]
The integers $c_{14}=D C_{14}(m/10^6)$ and $c_{12}=D C_{12}(m/10^6)$ satisfy the two angular tests precisely when
\begin{equation}\label{eq:integer-comparisons}
c_{14}H>DN,\qquad c_{12}H<DN,
\end{equation}
respectively. Here the letters $A,B,C$ denote endpoint numerators, not the angle sums of Section~\ref{sec:definitions}.

\subsection{Bounds for \texorpdfstring{$\pi$}{pi} and complete angle sums}
The only numerical bounds on $\pi$ used by the certificates are
\begin{equation}\label{eq:pi-bounds}
3.141592<\pi<3.141593.
\end{equation}
The replay verifies them with rational arithmetic. Machin's identity is
\[
\pi=16\arctan(1/5)-4\arctan(1/239).
\]
For $d=5,239$, the first six terms of the alternating series for $\arctan(1/d)$ give a strict lower bound, and adding the seventh gives a strict upper bound. Substitution proves \eqref{eq:pi-bounds} by rational comparisons. The identity follows from the tangent addition formula: $\tan(4\arctan(1/5))=120/119$, and subtracting $\arctan(1/239)$ gives tangent $1$ with the resulting angle in $(0,\pi/2)$.

Thus a sum of certified lower angle endpoints greater than $6283186$ microradians is greater than $2\pi$. A sum of certified upper endpoints less than $6283184$ microradians is less than $2\pi$. All angular margins reported in this paper are measured beyond these outer bounds, not beyond a floating-point value of $2\pi$.

\subsection{Independent verification of the exceptional certificate}\label{app:replay}
The exceptional-corona certificate is a depth-first trace of the interval proof. The replay program reconstructs the current box, checks the monotonicity direction of each contraction, certifies each supplied angle endpoint by \eqref{eq:integer-comparisons}, and checks the strict angle-sum inequality. It also checks that every split is internal, both children are visited, and no box remains at the end.

The proof trace cannot select another system: the replay hard-codes the seventeen lists in Table~\ref{tab:flowers}, the thirty-seven independent labels, and the initial integer endpoints $65105,100000$. The generator's sine-based trigonometric proposals are not trusted by the cosine-based replay. No floating-point operation or trigonometric library function is used by the latter.

\begin{center}
\begin{tabular}{@{}lr@{}}
\toprule
Verified quantity & Total\\
\midrule
Processed boxes & $7789$\\
Binary splits & $3894$\\
Rejected leaves & $3895$\\
Radius-bound contractions & $834176$\\
Strict angle-sum inequalities & $838071$\\
Pending or unresolved boxes & $0$\\
Minimum strict sum margin & $11\cdot10^{-6}$ radians\\
\bottomrule
\end{tabular}
\end{center}
The complete tree has $3895=3894+1$ leaves and $7789=2(3894)+1$ boxes, as also checked by the replay. Lemma~\ref{lem:certificate-soundness} applies to the verified trace and supplies the finite infeasibility statement in Proposition~\ref{prop:finite-exclusion}.

For example, the first step certifies
\[
T(84946,65105,65105)>897601/10^6.
\]
Since $7\cdot897601=6283207>6283186$, the equation at center $0$ forces its radius to exceed $0.84946$. The next step certifies
\[
T(70131,100000,100000)<1256631/10^6.
\]
Since $5\cdot1256631=6283155<6283184$, the equation at center $1$ forces its radius to be less than $0.70131$. Scale invariance explains the integer radius arguments in these examples. The full trace, not only these two entries, is needed for the proof.

\section{Reference constants and supplementary files}\label{app:supplement}
\subsection{Thresholds and exact reference comparisons}\label{app:reference}
The thresholds used by the proof are summarized below. They are rational inputs or rational enclosures and must not be substituted for one another.

\begin{center}
\small
\begin{tabular}{@{}p{.63\linewidth}p{.29\linewidth}@{}}
\toprule
Purpose & Value or interval\\
\midrule
Published degree restriction & $0.62$\\
Published degree-five adjacency and neighbor restrictions & $0.64$\\
Published two-face degree-seven exclusion & $0.645$\\
Finite exclusions and structure theorem: $q_0$ & $65105/100000$\\
Isolating interval for the exact optimum $s$ & $(0.65105018,0.65105019)$\\
Auxiliary, weaker bracket for $s$ & $(0.65104,0.65106)$\\
Bracket for the medium reference radius $r$ & $(0.83420,0.83440)$\\
\bottomrule
\end{tabular}
\end{center}
In particular, $0.65104$ is used only to weaken a lower bound on the fixed reference radius $s$ in scalar comparisons. It is not the cutoff used in the finite exclusions.

Exact evaluation of $P$ at $0.65105018$ and $0.65105019$ gives opposite signs. A positive rational lower bound for $P'$ on $[0.65104,0.65106]$ proves uniqueness of the root there. For this root, monotonicity and the exact endpoint comparisons
\[
T(s,1,1)+4T(s,1,0.83420)<2\pi,
\qquad
T(s,1,1)+4T(s,1,0.83440)>2\pi
\]
give the medium-radius bracket. These are checked conservatively by replacing the central radius $s$ by $0.65104$ in the upper comparison at the lower medium endpoint, and by $0.65106$ in the lower comparison at the upper medium endpoint.

For Table~\ref{tab:reference-bounds}, use a central medium radius at most $0.83440$, neighboring medium radii at least $0.83420$, and small radii at least $0.65104$. The sine enclosures of Lemma~\ref{lem:sine-enclosure} give the lower sums $6.333966$, $6.391892$, and $6.331199$. These establish all strict comparisons needed by Proposition~\ref{prop:reference-excess}. The exact flat reference flowers still come from the construction and angle identities in Section~\ref{sec:reference}, not from approximate closure in a drawing.

\subsection{Reproducing the finite proofs}\label{app:files}
The supplementary package contains the mathematical source, editable vector figures, the retained local programs, the complete exceptional-corona proof trace, and its separate checker. Their roles are as follows.
\begin{center}
\small
\begin{tabular}{@{}p{.40\linewidth}p{.54\linewidth}@{}}
\toprule
Program or data & Mathematical role\\
\midrule
\texttt{02\_pairing.js} & Proposition~\ref{prop:local}(i), including all finite degree cases.\\[3pt]
\texttt{03\_no\_three\_small.js} & Proposition~\ref{prop:local}(ii), using the pairing rule.\\[3pt]
\texttt{04\_reference\_checks.js} & Root and radius brackets, scalar angle comparisons, and the insertion inequality.\\[3pt]
\texttt{check\_patch\_combinatorics.py} & Face incidences, completion order, and source-disk counts for Table~\ref{tab:flowers}.\\[3pt]
\texttt{certificate.jsonl.gz} & Complete interval proof trace for \eqref{eq:finite-system}.\\[3pt]
\texttt{verify\_certificate.py} & Integer-only independent arithmetic replay.\\
\bottomrule
\end{tabular}
\end{center}
The local programs use JavaScript's arbitrary-precision integer arithmetic and require Node.js. The separate replay uses only the Python standard library. From the root of the supplementary package, all required checks can be run by
\begin{verbatim}
bash verify_all.sh
\end{verbatim}
The exceptional proof alone is verified by
\begin{verbatim}
cd certificates/exceptional
python3 check_patch_combinatorics.py
python3 verify_certificate.py certificate.jsonl.gz
\end{verbatim}
The geometry of the unfolded reduction is proved in Lemmas~\ref{lem:link} and~\ref{lem:geometric-reduction}; a successful combinatorial check does not replace that argument. Likewise, replay of the new certificate does not replace the published inputs in Proposition~\ref{prop:CZ} or the geometric coverage of the older local exclusions.

The retained local scripts are unchanged, so the auxiliary program also performs several comparisons not needed for the present proof. These additional outputs are not hypotheses of the theorems. The obsolete upper-bound search and the arguments about sparse exceptional sets are not used. The trust base of the computer-assisted part consists of the stated geometric reductions, the proved interval rules, the short implementations of those rules, and exact integer semantics. The programs are not a formal proof-assistant development.

\end{document}